
\input amstex
\documentstyle{amsppt}
\input xy
\xyoption{all}
\CompileMatrices
\xyReloadDrivers

\def\t#1 {\text{\rm #1}}
\def\b#1 {\text{\rm\bf #1}}


\topmatter
\title  Hyperelementary assembly for $K$-theory of virtually abelian groups
  \endtitle
  \rightheadtext{Hyperelementary assembly}
\author Frank Quinn\endauthor
\date September 2005, revised February 2006\enddate
\address Math, Virginia Tech, Blacksburg VA 24061-0123 USA\endaddress
\email quinn\@math.vt.edu\endemail
\thanks Partially supported by the National Science Foundation\endthanks
\abstract Controlled $K$-theory  is used to show that algebraic $K$-theory of virtually abelian groups is described by an assembly map defined using  possibly-infinite hyperelementary subgroups. The part coming from infinite subgroups is called the Farrell-Jones summand. This is shown to split from the finite-isotropy part; is parameterized by the rational projective space of the group; and a reduced version is shown to be torsion.
\endabstract

\endtopmatter

\head 1. Introduction\endhead
A group $G$ is $p$-{\it hyperelementary\/} for a prime $p$ if there is an exact sequence
$$1\to\text{ cyclic }\to G\to \text{ finite $p$-group }\to 1.$$
The finite ones are well-known and fundamental to the study of $K$-theory of finite groups. To these we add the possibility
that the cyclic group might be infinite. Infinite hyperelementary groups are simpler  than the finite ones: unless $G$ maps onto the infinite dihedral group it is a semidirect product $P\rtimes_{\alpha} T$ where  $T$ is infinite cyclic, $P$ is a $p$-group, and $\alpha$ is an automorphism of $P$ of order a power of $p$. Such $\alpha$ are nilpotent. Quite a number of classical methods for studying $K$-theory of finite groups (twisted group rings, for instance) should adapt to infinite hyperelementary groups.

The ``hyperelementary assembly conjecture'' suggests that the $K$-theory of group rings can be reconstructed from homology and the $K$-theory of hyperelementary subgroups using the assembly maps defined in \cite{Q1}. This is a sharpened version of the Farrell-Jones fibered isomorphism conjecture for $K$-theory, see \cite{FJ1, L, DL, BR}, and for an extensive survey \cite{LR}. It is also sharper (uses fewer subgroups) than an assembly result of Bartels and L\"uck \cite{BL}. The main result of this paper, Theorem 1.2.2, is that virtually abelian groups satisfy this conjecture. 

This result is new even in the simplest nontrivial cases  (product of a finite group and an infinite cyclic group), so it provides new calculations of Whitehead and $K_0$ groups. It also implies a version of the $K$-theory ``Novikov conjecture'' (5.1). However the main significance is the systematic nature of the description and its geometric implications. 

This paper completes a program first envisioned in \cite{Q2}, extending the pioneering work of Farrell-Hsiang \cite{FH1, FH2} in the torsion-free case.  The ingredient needed for \cite{Q2} to work is the fully-algebraic spectrum-valued controlled $K$-theory now provided by \cite{Q1}. Other versions of controlled $K$-theory have been developed by Pedersen-Weibel, Ferry-Pedersen and others, see the survey by Pedersen \cite{P}. Bartels and Reich \cite{BR} have a more recent  version. With the exception of \cite{BR} these are not effectively spectrum-valued, and all of them are missing the stability theorem vital to the present work. 

We call the cofiber of the assembly map from finite to virtually cyclic isotropy the ``Farrell-Jones component''. Theorem 1.2.2 shows this is the same as the cofiber from finite hyperelementary to all hyperelementary isotropy, and Theorem 1.3.3 shows it splits as a direct summand. Geometrically the finite-isotropy summand of $K$-theory is understood to relate to rigidity for manifold stratified models for classifying spaces, for instance $T\backslash L/\Gamma$ for $L$ a Lie group, $T$ a maximal torus, and $\Gamma$ a discrete subgroup, see \cite{Q3, CK}. This stratified rigidity is a generalization of the ``Borel conjecture'' that allows torsion in the group.  The Farrell-Jones component  seems to be related to codimension-two phenomena, but is still essentially mysterious.
In \S1.4 it is expressed as homology of a ``boundary'' for the group, with the discrete topology. This suggests a more coherent understanding might come from topologies, perhaps $p$-adic, on this boundary. 

Careful statements of the conjecture and theorem are given in \S\S1.1-1.2. The proof  is outlined in \S1.5 and given in \S\S2, 3. Group-theoretic lemmas are collected in \S4. 
Details for statements in \S1.4 about the Farrell-Jones component are given in \S5.  \S 6 provides a brief discussion of universal spaces.

I would like to thank Arthur Bartels, Jim Davis, Tom Farrell, Ian Hambleton, Wolfgang L\"uck, and the referee for pointing out errors in earlier versions of this paper. I am particularly grateful to the referee for help in clarifying many of the arguments.

\subhead 1.1 Assembly maps\endsubhead
In \cite{Q1} a controlled $K$-theory spectrum $\Bbb K(X;p,R)$ is defined for a map $p\:E\to X$ with $X$ a metric space and $R$ a commutative ring with 1. When the control space is a point this is equivalent to the standard algebraic $K$ spectrum of the group ring $R[\pi_1E]$.  When $p$ is a stratified system of fibrations the functor $\Bbb K$ can be applied fiberwise to get a stratified system of spectra over $X$, and there is an assembly map defined on a homology spectrum with coefficients in this system going to the controlled spectrum,
$$\Bbb H(X;\Bbb K(p;R))\to \Bbb K(X;p,R).\tag{1.1.1}$$
The Controlled Assembly Isomorphism Theorem of  \cite{Q1} asserts that this is an equivalence of spectra.  The result is formulated using {\it spectra\/} rather than {\it groups\/} because it is slightly stronger, more convenient, and notationally simpler. Note we are using compactly supported homology here, not the locally finite theory used in most of  \cite{Q1}. 

The uncontrolled assembly map
$$\Bbb H(X;\Bbb K(p;R))\to \Bbb K(\t pt ;E\to \t pt , R)\simeq \Bbb K(R[\pi_1(E,e_0)])\tag{1.1.2}$$
can now be described two ways. First, the homology construction is natural with respect to morphism of control maps, so we can apply it to  the morphism  
$$\xymatrix{E\ar[d]^{p}\ar[r]^=& E\ar[d]\\
X\ar[r]&{\t pt }.}$$
Second, if we identify the domain with controlled objects using 1.1.1 then we can relax the control from $X$ to  a point. In either case the map naturally goes to the spectrum $\Bbb K(\t pt ;E\to \t pt , R)$. 
Choosing a basepoint $e_0\in E$ gives an identification of this with the standard algebraic $K$-spectrum $\Bbb K(R[\pi_1(E,e_0)])$. For clarity we write it this way, though there are cases where basepoints cannot be chosen canonically enough to do this literally, see 5.1.3. This identification also requires $E$ to be connected, and this is not required when using the more natural notation. 

 In this paper we investigate cases where the uncontrolled assembly (1.1.2) is an isomorphism. The controlled assembly isomorphism theorem is used as a recognition criterion for this.

\subhead 1.2 Assembly conjectures\endsubhead
Suppose $\Cal G$ is a collection of isomorphism classes of groups. Given a group $G$ let $E_{G,\Cal G}$ denote a universal
(classifying) space for $G$-CW complexes with isotropy subgroups in $\Cal G$, see
\cite{FJ1, L} and \S6. Extreme cases are: 
\roster\item $\Cal G=\{1\}$ (trivial groups), then $E_{G, \{1\} }$ is the usual  free contractible $G$-space; and
\item if $G\in\Cal G$ then 
$E_{G, \Cal G }$ is equivariantly contractible.
\endroster
$G$ does not act freely on $E_{G,\Cal G}$ when $\Cal G$ is not trivial. The quotient is a stratified space stratified by orbit type: given a subgroup $H\in \Cal G$ the corresponding stratum in  $E_{G,\Cal G}/G$ is the image of points whose isotropy group  is exactly $H$. Note that acting by an element $h\in G$ changes the isotropy group by conjugation by $h$. The $G$-invariant subset is therefore points with isotropy conjugate to $H$. The stratum can be described two ways: $\{x\mid G_x \text{ conjugate to }H\}/G$, or $\{x\mid G_x=H\}/N(H)$ where $N(H)$ denotes the normalizer of $H$ in $G$. 

Geometric models for certain  $E_{G,\Cal G}$ will be important here. If $G$ is finitely generated virtually abelian then there is a projection with finite kernel to a crystallographic group. The crystallographic group has an action by isometries on $\b R ^n$. The induced action of $G$ on $\b R ^n$ is a universal space with finite isotropy groups, $ E_{G,\t finite }$. 

Now suppose $F$ is a free $G$-space. The quotient of the projection to the second factor gives a map 
$$\xymatrix{(F\times E_{G,\Cal G})\ar[r]^-{p_F}/G&E_{G,\Cal G}/G}.$$
This is a stratified system of fibrations over the orbit type stratification of the target. Note the fiber over a point with isotropy $H$ is $F/H$.
The discussion in 1.1 applies and associates assembly maps to this.  Assembly maps in this context are more commonly described using  Bredon homology c.f\. \cite{LR}. This is technically simpler, emphasizes the supportive context of group actions, and
Talbert \cite{T} has shown it gives the same theory.  We use the more general context for consistency with \cite{Q1} and to avoid having to translate in applications outside group actions.
\subsubhead 1.2.1 Definition\endsubsubhead We say that $K$ {\it is assembled from hyperelementary subgroups of\/}  $G$ if for every
free $G$-space $F$ and commutative ring $R$ the assembly map 
$$\Bbb H(E_{G,\t h.elem }/G;\Bbb K(p_F;R))\to \Bbb  K(R[\pi_1 F/G])$$
associated to the map $p_F\: (F\times E_{G,\t h.elem })/G\to E_{G,\t h.elem }/G$
is a weak equivalence of spectra. Here $E_{G,\t h.elem }$ is the universal space for actions with (possibly infinite) hyperelementary isotropy.

The free $G$-space $F$ is {\it not\/} assumed to be contractible or simply-connected, so this is the ``fibered'' form in the terminology of \cite{FJ1}. There is an exact sequence 
$$\xymatrix{1\ar[r]& \pi_1(F)\ar[r]&\pi_1(F/G)\ar[r]^-{\gamma} &G\ar[r]& 1}$$
Stalks in the coefficient system in the homology are $K$-theory of fibers in the stratified fibration, so equivalent to $\Bbb K(R[\gamma^{-1}(H)])$ for hyperelementary subgroups $H\subset G$. 
The conjecture therefore not only describes the $K$-theory of $G$ itself, but provides a reduction of the $K$-theory of any group mapping to $G$.  

The main result of the paper is:
\proclaim{1.2.2 Theorem}  $K$-theory of virtually abelian groups is assembled from hyperelementary subgroups.\endproclaim

The Farrell-Jones fibered isomorphism conjecture for $K$-theory, \cite{FJ1}, is the conjecture that $K$-theory is assembled from  virtually cyclic subgroups.  Hyperelementary is a subset of virtually cyclic so the hyperelementary conjecture is formally stronger. However virtually cyclic groups are included in Theorem 1.2.2 so the Farrell-Jones fibers can be further reduced (c.f.\ \S2.3) to their hyperelementary subgroups:
\proclaim{1.2.3 Corollary}  For a given $G$ the Farrell-Jones fibered isomorphism conjecture and  the hyperelementary assembly conjecture are equivalent.\endproclaim

In the finite case there is torsion information and this also extends to the infinite case.
If $\Lambda$ is a subring of the rationals then a $\Lambda$--{\it hyperelementary} group is either cyclic or $q$-hyperelementary for $q$ a prime not a unit in $\Lambda$. Note that cyclic groups are $q$-hyperelementary for all $q$, so explicitly including cyclic groups makes a difference only when $\Lambda=\b Q $ and there are no non-unit $q$. 
\proclaim{1.2.4 Theorem}    If\/ $G$ satisfies the hyperelementary assembly conjecture (i.e\. the $K$-theory of\/ $G$ is assembled from all hyperelementary subgroups), then $K\otimes \Lambda$ is assembled from $\Lambda$-hyperelementary subgroups. \endproclaim
This follows formally from the case with $G$ hyperelementary, and is proved in section 3.2.

\subhead 1.3 The Farrell-Jones component\endsubhead
\subsubhead 1.3.1 Definition\endsubsubhead The {\it Farrell-Jones component\/} $ \Bbb F\Bbb J(G;F,R)$ of virtually-cyclic isotropy homology is the cofiber (in the category of spectra) of the map from finite isotropy homology (or equivalently from finite hyperelementary to all hyperelementary):
$$\Bbb H(E_{G,\text{finite}}/G;\Bbb K(p_{\text{(finite)}},R))\to \Bbb H(E_{G,\text{v.cyc}}/G;\Bbb K(p_{\text{(v.cyc)}},R))\to \Bbb F\Bbb J(G;F,R).$$
Here, as above, $F$ is a free $G$-space, $p_{\text{(finite)}}$ is the projection $(E_{G,\text{finite}}\times F)/G \to E_{G,\text{finite}}/G$, and $\Bbb K(p_{\text{finite}},R)$ is the spectral cosheaf obtained by applying $\Bbb K$ to this fiberwise. 

In \S5  the iterated homology identity is used to rewrite  the Farrell-Jones component as homology with coefficients  the Farrell-Jones spectrum applied fiberwise:
\proclaim{1.3.2 Proposition} 
$$ \Bbb F\Bbb J(G;F,R)=\Bbb H(E_{G,\text{h.elem}}/G; \Bbb F\Bbb J(p_{\text{h.elem}},F,R))$$
\endproclaim 
This means that an understanding of the Farrell-Jones component in general will come from a functorial understanding of the Farrell-Jones component of infinite hyperelementary groups. The first application of this is the construction of a ``relaxation'' map that gives a left inverse for the map on homology in definition 1.3.1. This gives an extension of a result of  Bartels \cite{B}:
\proclaim{1.3.3 Proposition} The cofibration sequence in Definition 1.3.1 has a natural splitting, so 
the virtually-cyclic isotropy homology splits as the product of finite-isotropy homology and the Farrell-Jones component.
\endproclaim
Thus the Farrell-Jones component can also be described as the fiber of the relaxation map. This gives a formulation similar to the nilpotent-matrix description of the Bass Nil groups. A baby version of the argument gives a version of the Novikov conjecture for $K$-theory, in Section 5.1.

\subhead 1.4 Rational projective spaces\endsubhead
The homology description of the Farrell-Jones component suggests opportunities for interesting interactions between geometric group theory and $K$-theory. Fibers of the Farrell-Jones coefficient system in 1.3.2 are trivial over points in $E_{G,\text{h.elem}}/G$ that have finite isotropy. Deleting these points therefore does not change the homology but gives a much smaller base space. In particular subgroups that are not commensurable are in different components. Recall subgroups are {\it commensurable\/} if their intersection is finite index in each, and for infinite cyclic groups this is equivalent to nontrivial intersection. There are usually a lot of these classes, and we  try to organize them. 
\subsubhead 1.4.1 Definition (Rational projective spaces)\endsubsubhead A {\it rational line\/} in a group $G$ is a commensurability class of infinite cyclic subgroups. The set of all such lines is the {\it rational projective space\/} of the group and is denoted $\b QP (G)$. 

Adding a basepoint (notation ${\b QP }_+ (G)$) allows a description  as a quotient of $G$: torsion elements go to the basepoint and non-torsion elements go to the lines determined by the infinite cyclic subgroups they generate. ${\b QP }_+ (G)$ is natural with respect to homomorphisms;   homomorphisms with finite kernel induce injections and ones with finite-index image induce surjections. 
\subsubhead 1.4.2 Examples\endsubsubhead
\roster
\item If $G$ is a crystallographic group we can think of the free abelian subgroup as a lattice in $\b R ^n$. The rational lines in $G$ then correspond to real lines in $\b R ^n$ that intersect the lattice. This gives a description of $\b QP (G)$ as ``rational points''   in the real projective space. 
\item Elements $x,y\in G$ generate the same line if $x^n=y^m$ for some $m,n$, and they are not torsion. The group with presentation $\langle x,y\mid x^n=y^m\rangle$ therefore is a universal example. The center of this group is cyclic generated by $x^n$, and the quotient by the center is a free product of finite cyclic groups, $T_n*T_m$.
 Elements of $G$ generating the same line as $x$ are inverse images of torsion elements in the free product. 
 \item the rationals $\b Q $ form a single rational line.
 \endroster
 In example (2) the set of elements generating subgroups commensurable to the center is not a subgroup (the torsion elements in the free product are not a subgroup). However we see that this phenomenon is connected to  free abelian subgroups of rank 2 (inverse of an infinite cyclic subgroup of the free product).

\subsubhead 1.4.3 Decomposing the Farrell-Jones component\endsubsubhead
The observation at the beginning of Section 1.4 can now be formalized: $G$ acts by conjugation on the projective space, and the Farrell-Jones component breaks into a sum over the conjugacy classes. Inspired by the crystallographic example we ask: are there metrics on the projective space so that the Farrell-Jones coefficient system (or parts of it) extend continuously over the completion? To see what this would entail we write a direct sum  as homology of the index set thought of as a space with the discrete topology. 
$$\Bbb F\Bbb J(G;F,R) = \Bbb H(\b QP (G)/G;\Bbb F\Bbb J(r;F,R)).$$
Here $\Bbb F\Bbb J(r;F,R)$  is a ``cosheaf'' whose stalk over a rational line is the corresponding summand in homology. This is a bit less straightforward than one might have hoped because these rational lines may not be subgroups. Extension of the coefficient system over a completion would give a  continuous version more likely to have a systematic description. In applications we could either restrict to the rational points (and get involved in trying to enumerate them) or we may be led to new geometric problems that use the topology.

The final remark is that a principal technique in algebra is to pass to the $p$-completion of the coefficient ring. There are hints that this will induce some sort of $p$-adic topology on the rational projective space. It seems possible that for large classes of groups these $p$-adic spaces would be algebraic varieties. This would be a very useful way to see systematic and functorial structure in $K$-theory. Description of commensurability classes as rational points in a variety is also an interesting variation on solvability of the word problem.

\subhead 1.5 Outline of the proof\endsubhead
A consequence of the Swan-Lam-Dress  induction theory for finite groups
\cite{D} is that  finite groups satisfy hyperelementary assembly. 
The main theorem is proved in Section 2 using induction in finite quotients and the controlled theory, a method developed by Farrell-Hsiang  \cite{FH1, FH2} for torsion-free poly-(finite or cyclic) groups.  

By general principles the theorem reduces to finitely generated groups. Next we use the fibered reduction feature. If $G$ has a normal infinite cyclic subgroup we can divide and reduce the rank. We can also pass to $G/H$ if $H$ is a finite normal subgroup. This reduces the theorem to some special cases discussed in the next paragraph, and crystallographic groups without normal infinite cyclic subgroups. These last have nice non-hyperelementary finite quotients. We can arrange the quotients so that hyperelementary subgroups are either smaller in group-theoretic terms or lie in the image of an expanding map with arbitrarily large preassigned expansion coefficient. In the first case we proceed by induction on group-theoretic size. In the second case the pullback of $K$-theory objects become very small in a metric sense, and the controlled theory identifies them as coming from homology via the assembly map. In either case pullbacks to hyperelementary subgroups come from assembly, so finite induction theory asserts that the whole $K$-theory comes from assembly. 

The reduction described above requires  rank--1 groups and a few rank--2 crystallographic groups to be done separately. These are treated in \S3. The approach is similar to the general case: find finite quotients so hyperelementary subgroups are either smaller in some group-theoretic sense or factor through an expanding map so the control theory applies.

\head 2. The main argument \endhead
In this section the proof that virtually abelian groups satisfy hyperelementary assembly is reduced to rank one and a few rank two cases. These cases are done in \S3. The proof is outlined in \S2.1. Basic reductions are given in \S2.2--2.4.

\subhead 2.1 Structure of the proof \endsubhead
After a reduction to finite generation in 2.2 the proof  centers on hyperelementary induction in finite quotients following a plan developed by Farrell and Hsiang \cite{FH1}. Suppose $1\to A\to G\to Q\to 1$ is exact, with $Q$ finite and $A$ finitely generated free abelian. Let $p$ be a prime not dividing $|Q|$ and $r$ a positive integer, then $G/A^{p^r}$ is usually not hyperelementary. Let $H\subset G/A^{p^r}$ be a hyperelementary subgroup, and $\hat H$ the preimage in $G$. According to the fibered reduction of 2.3 if the conjecture is true for all such $\hat H$ then it is true for $G$.  There are three cases: 
\roster\item the projection $H\to Q$ is not onto;
\item $H\to Q$ is onto but not injective; or
\item $H\to Q$ is an isomorphism.
\endroster
We want to avoid case (2)  in the main argument. 
If $p$ is chosen carefully then 2.4.2 shows that in case (2) there is an infinite cyclic normal subgroup of $G$, i.e.\ $G$ is cyclically reducible. Dividing gives a map to a group of lower rank. If we proceed by induction on rank then the result is known for the quotient. It then follows for $G$ by the fibered reduction of 2.3 and a few rank-2 special cases that must be done separately, see \S2.4. Therefore  we may assume $G$ is not cyclically reducible and case (2) does not occur. 

In case (1) $\hat H$ has smaller finite quotient than $G$.  If we proceed by induction on the order of the finite quotient then we can assume the result is known for these $\hat H$. 

Case (3) is where geometry is used. A minor reduction in 2.3 shows we can suppose $G$ to be a crystallographic group. The standard action of $G$ on $\b R ^n$ is a model for the universal $G$-space with finite isotropy. It is enough to show that the assembly map associated to $(F\times \b R ^n)/G\to (\b R ^n)/G$ is an equivalence since the hyperelementary version follows from this and hyperelementary assembly for finite groups. Use the controlled version \cite{Q1} to identify the domain of the finite-isotropy assembly as the subset of the $K$ space $\Bbb K(R[\pi_1(F\times \b R ^n/G)])$ that is $\epsilon$ controlled over the quotient $\b R ^n/G$ for sufficiently small $\epsilon$. The objective is to show that the whole space deformation retracts to the controlled subspace. The mechanism for this is transfers by expansive maps. 

Suppose $p^r$ is congruent to 1 mod $|Q|$. According to 2.5 there is an expansive homomorphism $\gamma\:G\to G$ with expansion $p^r$ and the preimage of a hyperelementary subgroup  $H\subset G/A^{p^r}$ factors through this up to conjugacy  in case (3) ($H\to Q$ an isomorphism). Multiplication by $p^r$ $\b R ^n\to \b R ^n$ is $\gamma$-equivariant. Transfer to $\hat H$ corresponds to pulling back by this map, so it {\it reduces\/} arc length in $\b R ^n$ and $\b R ^n/G$ by a factor of $p^r$. If we fix an initial size then using a really big power of $p$ gives result with arc length less than $\epsilon$ and therefore controlled. This is enough to show the assembly is a weak equivalence. Weak equivalence requires that something happen for finite complexes mapping into the space. A map of a finite complex into a $K$-space is determined by a finite number of paths and homotopies, so there is an upper bound on the length and some fixed power of $p$ will shrink it below the stability threshold.  

This completes the basic outline. We describe a modification that shows the assembly for finitely generated $G$ is an equivalence, not just a weak equivalence. This also gives a bit more detail for the final step just above.

To show the assembly is a genuine equivalence we use induction on another measure of size. Simplices of the $K$-space are defined in terms of paths and homotopies in $(F\times \b R ^n)/G$. By small approximation we may assume these are smooth in the $\b R ^n$ coordinate.  Let $\Bbb K_{\ell}$ for $\ell>0$ denote the subspace of simplices in which all the paths (including the ones in the homotopies) have arc length less than $\ell$. We suppose as an induction hypothesis that there are $\ell_{k-1}>k-1$ and $\ell_k>k$ and a homotopy of $\Bbb K_{\ell_{k-1}}$ rel $\Bbb K_{\epsilon}$ into $\Bbb K_{\epsilon}$, and the homotopy stays in $\Bbb K_{\ell_k}$.

 We are ignoring some routine $\epsilon-\delta$ details here: ``$\epsilon$'' generically represents numbers so small that the stability theorem applies for $\epsilon$ and $\b R ^n/G$. Also the stability theorem uses ``radius'' rather than arc length as the measure of size. Radius is smaller so it is sufficient to bound length, and length is more appropriate for large objects. For instance any arc mapping to a circle has radius $<1$, but it may have arbitrarily large arc length and this is what is significant when pulling back to covers. 
 
 The objective is to find $\ell_{k+1}>k+1$ and an extension of the homotopy to $\Bbb K_{\ell_k}$. Choose a good $p$ as in 2.4 and $r$ so that $p^r\equiv 1$ mod $|F|$ (to get an expansive map) and so large that $\ell_k/p^r <\epsilon$. Dress induction implies that the identity map of $\Bbb K$ for $G$ is a linear combination of inclusions of transfers to $\hat H$, for inverses of hyperelementary 
$H\subset G/A^{p^r}$. It is therefore sufficient to show that the transfers have extensions of the homotopy. If $H$ is in case (1) then the conjecture is assumed known as an induction hypothesis and the homotopy extends to the whole space. Since $G$ is cyclically irreducible case (2) does not occur. Finally in case (3) $\hat H$ is the image of an expansive homomorphism with expansion $p^r$. We have chosen $p^r$ so that the whole subspace $\Bbb K_{\ell_k}$ pulls back into  $\Bbb K_{\epsilon}$ and the stability theorem gives the desired extension of the homotopy.

Some care is required to see that homotopies obtained this way are contained in some  $\Bbb K_{\ell_{k+1}}$, so the induction can be continued. For the part coming from the expansive construction this is part of the stability theorem. The other part can be handled rather crudely since there are (up to isomorphism) only finitely many groups with smaller finite quotient. We can take the max of estimates for all such groups.

\subhead 2.2 Finite generation \endsubhead
The following is standard and formal:
\proclaim{2.2.1 Proposition} The assembly map of 1.1.4 is weakly equivalent to the (homotopy) direct limit of assembly maps for finitely generated subgroups of $G$, or finitely presented groups mapping to $G$. \endproclaim
\subsubhead Reason\endsubsubhead
Weak equivalence means they behave the same with respect to maps of finite complexes. On the left side, a simplex in $\Bbb K(\t pt ;E\to \t pt, R)$  is determined by a finite number of paths and maps of 2-cells in $E$. A map of a finite complex is also determined by such finite data, so lies in a finite subcomplex of $E$. This has finitely generated fundamental group. The argument for the right side is similar if we use the description of homology as controlled $K$-theory. Abbreviate the control map $(F\times E_{G,\t he })/G\to E_{G,\t he }/G$ as $E\to X$, then a simplex of the space is an inverse limit of $K$-objects in $E$ with $\epsilon$ control in $X$ with $\epsilon\to 0$. The stability theorem asserts that the limit is determined by the value at some finite $\epsilon$, and this is described by a finite number of paths and homotopies. Again it is supported in finite subcomplexes of $E$ and $X$.
\proclaim{2.2.2 Corollary} Hyperelementary assembly for all virtually abelian groups follows from  the finitely generated case.\endproclaim
The finite-index abelian subgroup of such a group is finitely generated. Note finiteness of the quotient is needed for finite generation of the subgroup: $\b Z [\frac12]\rtimes  T $ where the generator of $T$ acts by multiplication by 2 is finitely generated, but the normal abelian subgroup is not.

\subhead 2.3 Fibered reductions \endsubhead
The beauty of the ``fibered'' formulation of the conjecture is that it satisfies an extension property, see Farrell-Jones \cite{FJ1, \S1.7}.
\proclaim{2.3.1 Proposition} Suppose $\gamma\: G_1\to G_2$ is a homomorphism, $G_2$ satisfies hyperelementary assembly, and for every hyperelementary subgroup $H\subset G_2$ the inverse image $\gamma^{-1}(H)$ satisfies hyperelementary assembly. Then $G_1$ satisfies hyperelementary assembly.\endproclaim
\subsubhead Proof\endsubsubhead
Suppose $F$ is a free $G_1$ space, and recall $E_{G,\t he }$ denotes the universal $G$-space with hyperelementary isotropy. Denote the sequence of spaces 
$$(F\times E_{G_1,\t he }\times E_{G_2,\t he })/G_1\to (E_{G_1,\t he }\times E_{G_2,\t he })/G_1\to (E_{G_2,\t he })/G_2\tag{2.3.2}$$
by $$\xymatrix{E\ar[r]^{p_1}&X_1\ar[r]^{g}&X_2}$$
and denote $gp_1$ by $p_2$. The maps in sequence 2.3.2 are projections on factors, and $G_1$ acts on $E_{G_2,\t he }$ via the homomorphism $\gamma$. The maps are nicely stratified so we get a commutative diagram 
$$\xymatrix{{\Bbb H}(X_1;\Bbb K(p_1;R))\ar[r]^{A_1}\ar[d]^{\simeq}&{\Bbb K}(R[\pi_1E])\\
{\Bbb H}(X_2;\Bbb H(g;\Bbb K(p_1;R)))\ar[r]&{\Bbb H}(X_2;\Bbb K(p_2;R))\ar[u]^{A_2}
}\tag{2.3.3}$$
Here $A_1$ and $A_2$ are assembly maps associated to $p_1$ and $p_2$ respectively, the left vertical map is the ``iterated homology identity'' of Part I \S6.10, and the bottom map is induced on homology by a morphism of coefficient systems over $X_2$. The iterated homology identity uses the coefficient system obtained by taking homology of point inverses of the map $g$ with coefficients in the restriction of the system over $X_1$. We will identify the maps in the diagram in terms of the isomorphism conjecture. 

First, the assembly $A_1$ is the one associated with $F$ that we want to show is a weak equivalence. The map $p_1$ in 2.3.2 differs from the standard formulation by having an additional factor. But projection off this factor, 
$$(E_{G_1,\t he }\times E_{G_2,\t he })/G_1\to (E_{G_1,\t he })/G_1$$
is stratum-preserving and fibers are $E_{G_2,\t he }$ and therefore contractible. The extra factor therefore does not change the $K$-theory or homology. 

Next $A_2$ is the assembly map associated to the free $G_2$ space 
$$G_2\times_{G_1}F=(G_2\times F)/(h,gx)\sim (h\gamma(g)^{-1},x).$$
Note that $(G_2\times_{G_1}F)/G_2=F/G_1$. By hypothesis $G_2$ satisfies hyperelementary assembly so $A_2$ is a weak equivalence.

The final task is to show that the lower map in diagram 2.3.3 is a weak equivalence. This is induced on homology by a morphism of coefficient systems over $X_2$ so it is sufficient to show this morphism is a weak equivalence over each point in $X_2$. Over a point $x$ the morphism is an assembly map
$$\Bbb H(g^{-1}(x);\Bbb K(p_1;R))\to \Bbb K(x; p_2^{-1}(x)\to x, R)\simeq \Bbb K(R[\pi_1 (p_2^{-1}(x))]).$$
In fact it is an assembly map associated to a subgroup of $G_1$. Recall that $g$ is defined by dividing the projection $E_{G_1,\t he }\times E_{G_2,\t he }\to E_{G_2,\t he }$
 by group actions. $g^{-1}(x)$ is therefore $E_{G_1,\t he }/\gamma^{-1}(H)$ , where $H$ is the isotropy group of a preimage of $x$ in $E_{G_2,\t he }$. By definition $H$ is hyperelementary, and it is a simple standard fact that $E_{G_1,\t he }$ is equivalent to $E_{\gamma^{-1}(H),\t he }$ as  a $\gamma^{-1}(H)$-space. Finally the inverse in $E$ is equivalent to $(F\times E_{\gamma^{-1}(H),\t he })/\gamma^{-1}(H)$, where $F$ is thought of as a free $\gamma^{-1}(H)$-space. The hypothesis that inverses in $G_1$ of hyperelementary subgroups in $G_2$ satisfy assembly therefore implies that this assembly map is a weak equivalence. 
 
This shows that the vertical and lower maps in diagram 2.3.3 are weak equivalences, so the upper map is also, as was to be proved.

An easy corollary gives the connection to geometry:
\proclaim{2.3.4 Corollary} If crystallographic groups and virtually cyclic groups satisfy hyperelementary assembly  then so do all finitely generated virtually abelian groups.\endproclaim
\subsubhead Proof\endsubsubhead Recall (\S4.2) that a virtually abelian group has a crystallographic quotient of the same rank, say $G\to \Gamma$.  The  kernel is finite so the inverse images of hyperelementary subgroups of $\Gamma$ are all virtually cyclic. According to the fibered reduction theorem 2.3.1 if $\Gamma$ and all these inverse images satisfy assembly then so does $G$. 

\subhead 2.4 Normal cyclic subgroups \endsubhead
We say a virtually abelian group is {\it cyclically reducible\/} if there is a normal infinite cyclic subgroup. In this case there is an exact sequence 
$$1\to T\to G \to \hat G\to 1$$
and the fibered reduction theorem gives assembly for $G$ from assembly for the lower-rank group $\hat G$ and inverse images. More precisely:
\proclaim{2.4.1 Proposition}Hyperelementary assembly  for virtually abelian groups follows from assembly for:
\roster\item rank-1 (i.e.\ virtually infinite cyclic) groups;
\item rank-2 crystallographic groups with a normal infinite cyclic subgroup; and
\item cyclically irreducible crystallographic groups.
\endroster\endproclaim
The main body of the proof, described in \S2.1, shows that the third case follows from the first two. 
\subsubhead Proof\endsubsubhead Inverses of finite hyperelementary subgroups of $\hat G$ are virtually cyclic. Inverses of infinite hyperelementary subgroups are rank two, and these can be reduced to  crystallographic groups as in 2.3.2. Thus assembly for $G$ follows from the groups in cases (1) and (2) and assembly for $\hat G$. If $\hat G$ is cyclically reducible then repeat, and continue until an irreducible quotient is reached. This happens after finitely many steps because the rank decreases in each step. Finally reduce to the crystallographic quotient as in 2.3.4. 

Much of the utility of this reduction comes from being able to recognize cyclic reducibility using finite quotients. Suppose 
$$1\to A \to G\to Q\to 1$$ 
is a virtually abelian presentation of $G$. The following is implicit in Farrell-Hsiang \cite{FH1}:
\proclaim{2.4.2 Proposition} Suppose $A$ is torsion-free and $p$ is a prime not congruent to 1 mod any odd prime divisor of\/ $|Q|$, and congruent to 3 mod 4 if\/ $|Q|$ is even. If there is a nontrivial\/ $Q$-invariant cyclic subgroup of\/ $A/A^{p^r}$ then there is an infinite cyclic $Q$-invariant cyclic subgroup of\/ $A$. Further if\/ $H\subset G/A^{p^r}$ is hyperelementary and $H\to Q$ is onto but not injective then the kernel contains such a nontrivial $Q$-invariant cyclic subgroup. \endproclaim
$Q$ acts on $A$ by conjugation in $G$, so a subgroup of $A$ is $Q$-invariant if and only if it is normal as a subgroup of $G$. The conclusion of the Proposition is therefore that $G$ is cyclically reducible.

Suppose $C\subset A/A^{p^r}$ is an invariant cyclic subgroup. The cyclic subgroup of $C$ of order $p$ is also invariant, so there is no loss in assuming $C$ has order $p$. 
\proclaim{2.4.3 Lemma} There is a homomorphism $\rho\:Q\to\{\pm1\}$ so that $Q$ acts on $C$ by $f\,  c=c^{\rho(f)}$.\endproclaim
Another way to put this is that $\{\pm1\}$ is the only $|Q|$-torsion in the automorphism group of $C$. This is well-known, but we give a short proof. 

\subsubhead Proof of the Lemma\endsubsubhead First suppose $Q$ is a $q$-group for $q$ an odd prime. Then orbits of the action of $Q$ on $C$ have order powers of $q$. The number of fixed points (orbits with $q^0$ elements) is therefore congruent to the order of $C$ mod $q$. Since we assumed $p$ is not congruent to 1 mod $q$ there must be a fixed non-1 element. But fixed elements form a subgroup and $C$ has prime order so $Q$ fixes all of $C$.

In general the argument above shows the Sylow subgroups of $Q$ of odd order fix $C$. The image of $Q$ in the automorphisms of $C$ must therefore be a 2-group. Assume $Q$ is a 2-group. Orders of orbits are powers of 2, so the order of $C$ mod 4 is congruent to the number of fixed elements plus the number of elements with orbits of order 2. Since we have supposed $p$ is congruent to 3 mod 4, either there are at least two nontrivial fixed elements or at least one element, say $c$, with orbit of order 2. In the first case $C$ must be fixed by $Q$. In the second case $Q$ acts on the element $c$ by $f\, c=c^{\rho(f)}$ for some homomorphism $\rho\:Q\to \{\pm1\}$. But the elements of $C$ on which $Q$ acts  by $\rho$ is a subgroup, and so $Q$ acts on all of $C$ in this way. This completes the proof of the Lemma.

\subsubhead Proof of 2.4.2\endsubsubhead
 From the lemma we know there is an element $a\in  A$ so that $Q$ acts on the image $[a]\in A/A^{p^r}$ by $[f\,a]=[a]^{\rho(f)}$. Define a new element by
 $$\hat a =\sum_{f\in Q}f\,a^{\rho(f)}.$$
Then\roster\item for any $g\in Q$, $g\, \hat a =a^{\rho(g)}$; and
\item the image $[\hat a]\in A/A^{p^r}$ is $a^{|Q|}$, so is nontrivial since $p$ is prime to $|Q|$.\endroster
This element generates the invariant infinite cyclic subgroup we seek. Note this is contained in a minimal summand of $A$ that is also $Q$-invariant.

For the final part of the proposition suppose $H$ is a $q$-hyperelementary subgroup of $G/A^{p^r}$ that maps onto $Q$. If $p\neq q$ then the kernel lies in the cyclic part of $H$ so is a cyclic normal subgroup of $H$. Since $H\to Q$ is onto, this subgroup is fixed under the conjugation action of $Q$ and therefore is also normal in $G/A^{p^r}$. If $p=q$ then $Q$ is the cyclic part of $H$ and $H=K\times Q$. In this case the whole kernel is invariant under the $Q$ action, so any cyclic subgroup is normal. 

\subhead 2.5 Finite groups \endsubhead
The following is an application of Dress' induction theory \cite{D1}, see Bartels-L\"uck \cite{BL, Lemma 4.1}.
\proclaim{2.5.1 Theorem} Finite groups satisfy hyperelementary assembly. More generally if $\Lambda$ is a subring of the rationals then for finite $G$ then  $K\otimes \Lambda$ is assembled from the $\Lambda$-hyperelementary subgroups of $G$.\endproclaim
See Theorem 1.2.4 for $\Lambda$-hyperelementary groups.

In its raw form Dress' theorem describes the homotopy groups of the $K$-theory spectrum as the direct limit of $K$-theory groups associated to hyperelementary subgroups of $G$. One then recognizes this direct limit as the homology group appearing as the domain of the assembly, and the direct limit map is the assembly. 
Dress actually gives a sharper result: The identity map of the  $K$ spectrum to itself is a linear combination of inclusions of transfers to hyperelementary subgroups. We will use this version to prove things  using geometric properties of transfers.

\subhead 2.5 Expansive maps \endsubhead
Crystallographic groups act as isometries on Euclidean spaces of the same rank, with compact quotient. Let $1\to A\to G \to Q\to 1$ have $A$ finitely generated free abelian and $Q$ finite acting faithfully on $A$. The splitting group $A\rtimes Q$ acts in the evident  way on  $A\otimes {\b R }\simeq {\b R }^n$, and $G$ embeds in the splitting group. If $r\equiv 1$ mod $|Q|$ and $\gamma\:G\to G$ is the expansive endormorphism with expansion $r$ as in  4.3, then multiplication by $r$ is a $\gamma$-equivariant map we denote by  $g\:{\b R }^n\to {\b R }^n$. 

The expansive map on ${\b R }^n$ induces a map of quotients ${\b R }^n/G\to {\b R }^n/G$. This is usually  a ``branched cover'' of stratified spaces. For instance if $G$ is the infinite dihedral group then ${\b R }^n/G$ is the unit interval $[0,1]$. The expansive map with odd expansion $r$ is $\xymatrix{[0,1]\ar[r]^r&[0,r]}$ followed by the folding map $[0,r]\to [0,1]$. When $F$ is a free $G$-space the induced map $(F\times {\b R }^n)/G\to (F\times {\b R }^n)/G$ is a covering space. 

The crystallographic action is a model for the universal space with finite isotropy.   Simplices in $\Bbb K(R[\pi_1(F/G)])$ are defined using paths and homotopies in $F/G\simeq (F\times {\b R }^n)/G$. Transfers by $\gamma$ are defined geometrically by pulling back these data over the quotient of the expansive map. Naturally we may approximate paths and homotopies by ones that are smooth in the ${\b R }^n$ coordinate, and in this case arc lengths are defined. The following is then clear:
\proclaim{2.5.1 Lemma} Pulling smooth paths back by an expansive map with expansion $r$ divides arc length in ${\b R }^n/G$ by $r$.\endproclaim
In particular if $r$ is really big compared to the length of an arc then the pullback is very short.
  
\head 3. Special cases \endhead
In \S2 the assembly theorem for virtually abelian groups is reduced to the rank 1 case, and a few rank 2 groups.  These special cases are treated here. The tools are the same as the main proof.

\subhead 3.1 Rank one \endsubhead
\proclaim{3.1.1 Proposition} Virtually cyclic groups satisfy hyperelementary assembly.\endproclaim
The proof is similar in outline to the main proof in 2.1. There is a sequence  $$1\to A\to G\to Q\to 1$$ with $A$ infinite cyclic and $Q$ finite. We proceed by induction on the order of $Q$, so suppose the result is known for groups with smaller quotient. After preliminary reductions we find a situation in which either something gets smaller or we can use an expansive map. 

If the order of $Q$ has a single prime divisor then $G$ is infinite hyperelementary and there is nothing to prove. Thus suppose $p$,  $q$ are distinct primes dividing the order of $Q$, $r$ is very large, and consider a hyperelementary subgroup  $H\subset G/A^{(pq)^r}$. According to 2.3.1 and 2.5.1 we need to show that the preimage $\hat H\subset G$ satisfies hyperelementary assembly.

We suppose $H\to Q$ is onto, since otherwise the preimage $\hat H\subset G$ has smaller finite quotient than $G$ and the result is already known. Finally $H$ is hyperelementary with respect to some prime so one of $p, q$ must be different from that prime. Suppose it is $p$, then the $p$-torsion of $H$ is a normal cyclic subgroup we denote as $H_p$. Note this implies $Q_p$ is also cyclic. 

Let $G_p$ denote a maximal $p$-torsion subgroup of $G$. The proof divides into three cases:\roster\item $Q_p$ acts on $A$; 
\item $G_p\to Q_p$ is onto and $Q_p$ does not act on $A$; and
\item $G_p\to Q_p$ is not onto.
\endroster
\subsubhead 3.1.2 $Q_p$ acts on $A$\endsubsubhead
In this case $p=2$ and $G$ has dihedral quotient $j\colon G\to D=T\rtimes T_2$. Note $j(A)=(T)^n$ for some $n$ so there is an induced map $\bar j\colon  G/A^{(2q)^r}\to \bar D= T/T^{n(2q)^r}\rtimes T^2$. 
The image $\bar j(H)$ is both $p$--hyperelementary (because $H$ is) and 2--hyperelementary (because $\bar D$ is) so it must be cyclic. But since the composite $\bar j(H)\to T_2$ is surjective this composite must be an isomorphism. It follows the the index of $j(H)\subset D$ is divisible by the index of $\bar j(H)\subset \bar D$ and therefore by $q^r$. This is an odd number so $j(H)$ is conjugate into the image of the expansive map of index $q^r$ on $D$, and pulling back to $H$ reduces lengths. By choosing sufficiently large $r$ we can reduce any given finitely generated part of the $K$-theory to finite subgroups, and consequently get the desired reduction for $\hat H$ and thus $G$. 

\subsubhead 3.1.3 $G_p\to Q_p$ an isomorphism\endsubsubhead
Let $\hat Q_p$ denote the preimage in $G$, then the exact sequence for $G$ gives 
$$1\to A\to \hat Q_p\to Q_p\to 1.$$
The subgroup $G_p\subset \hat Q_p$ gives a splitting of this sequence, and by 3.1.2 we may assume $Q_p$ acts trivially on $A$, so $\hat Q_p=A\times G_p$ (direct product of subgroups). Further $\hat Q_p/ A^{(pq)^r}=(A/A^{(pq)^r})\times \pi(G_p)$ where $\pi\colon G\to G/A^{(pq)^r}$ is the natural epimorphism. 

Recall we are considering an  $\ell$-hyperelementary subgroup $H\subset G/A^{(pq)^r}$ such that $H\to Q$ is onto and $p\neq \ell$. The $p$-torsion in $H$ is the intersection with $\hat Q_p/A^{p^r} =A /A^{p^r}\times \pi(G_p)$.  The $p$-torsion in $G$ injects to $\{1\}\times Q_p$, so if $H\cap\pi(G_p)$ is trivial then $\hat H$ has no $p$-torsion. In this case $\hat H\to Q/Q_p$ has infinite cyclic kernel so $\hat H$ has smaller quotient than $G$ and the result is already known. The problematic cases are therefore those in which $H\cap \pi(G_p)\neq \{1\}$.

Choose a generator of the normal cyclic subgroup prime to $\ell$. Express the image in $A /A^{p^r}\times \pi(G_p)$ of this generator as $a^mb^n$ where $a,b$ are generators of $A$, $\pi(G_p)$ respectively. Denote the power of $p$ dividing the order of $H$ by $p^s$, then there is a nontrivial element of $H\cap \pi(G_p)$ provided there is $j$ so that $a^{mj}=1$ and $b^{nj}\neq 1$. This is equivalent to requiring that $p^r$ divides $mj$ and $p^s$ does not divide $nj$. This implies that $p^{r-s}$ divides $m$, and consequently $H$ lies in $A^{p^{r-s}}/A^{p^r}\times Q_p$. This in turn implies that the image of $\hat H$ in the crystallographic quotient of $G$ lies in the image of the expansive map of index $p^{r-s}$. Recall that $s$ is fixed (the exponent of the order of $Q_p$) and $r$ is arbitrarily large, so as usual we can use this to reduce $K$-theory to finite subgroups and obtain the result. 

\subsubhead 3.1.4 $G_p\to Q_p$ not an isomorphism\endsubsubhead
This is the final case in the proof of 3.1.1, and is ``soft'' in that it does not require expansive maps or control. We will see that the preimage $\hat H$ of the hyperelementary subgroup is disjoint from the $p$-torsion in $G$, so $\hat H\to Q/Q_p$ also has  cyclic kernel and has quotient smaller than $Q$. 

As in 3.1.3 we begin with the sequence  $1\to A\to \hat Q_p\to Q_p\to 1$. $Q_p$ acts trivially on $A$ so the crystallographic quotient is infinite cyclic: $\hat Q_p\to T$. Choosing a splitting $\hat T\subset \hat Q_p$ gives a decomposition $A=G_p\rtimes_{\alpha}\hat T$ where $\alpha$ is an automorphism of $G_p$. 

Note that dividing $\hat Q_p$ by both $A$ and $G_p$ gives an identification $T/A=Q_p/G_p$. Denote the orders of $G_p$ and $Q_p/G_p$ by $p^s$, $p^t$ respectively, and note that we are assuming $s$, $t$ are both greater than 0. It follows that ${\hat T}^{p^{s+t}}=A^{p^s}$, so in particular if $r\geq s$ then $\hat Q_p/A^{p^r}=G_p\rtimes \hat T/\hat T^{p^{r+t}}$. 

Now return to the $\ell$-hyperelementary group $H\subset G/A^{(pq)^r}$, with $\ell\neq p$. The $p$--torsion of $H$ is the cyclic group $H_p\subset \hat Q_p/A^{p^r}=G_p\rtimes \hat T/\hat T^{p^{r+t}}$. Recall that the goal is to show $\hat H$ is disjoint from the $p$--torsion of $G$, or equivalently that $H_p$ intersects trivially the image of the $p$-torsion, $\{1\}\rtimes G_p$. This image is the kernel of the projection $G_p\rtimes \hat T/\hat T^{p^{r+t}}\to T/T^{p^{r+t}}$ so we want to show $H_p\to T/T^{p^{r+t}}$ is injective. Recall that we assume $H\to Q$ is onto, so $H_p\to Q_p/G_p=T/T^{p^{t}}$ is onto. But this is a nontrivial quotient of the cyclic $p$--group $T/T^{p^{r+t}}$, so $H_p\to T/T^{p^{r+t}}$ is also onto. The injectivity we seek is therefore equivalent to showing that $H_p$ has order at most ${p^{r+t}}$. This is done by showing that the group $G_p\rtimes T/T^{p^{r+t}}$ has exponent ${p^{r+t}}$.

Recall that ${\hat T}^{p^{s+t}}=A^{p^s}$ and this commutes with $G_p$. The endomorphism $\alpha$ in the semi-direct product  is defined by conjugating by a generator of $\hat T$ so it has order dividing $p^{s+t}$. Since $G_p$ is a cyclic $p$-group this implies $\alpha$ is of the form $\alpha(x)=x^{1+pn}$ for some integer $n$. Now suppose $xt^j\in G_p\rtimes T/T^{p^{r+t}}$  and compute
$$(xt^j)^p=xt^jxt^j\dots =(\prod_{i=0}^{p-1}\alpha^{ij}(x))t^{pj}=x^wt^{pj}$$
where $w=\sum_{i=0}^{p-1}(1+pn)^{ij}$. Each term $(1+pn)^{ij}$ has the form $1+pm$ for some $m$, and there are $p$ terms in the sum so the sum is divisible by $p$. Thus we see $(xt^j)^p=x^{pm}t^{pj}$ for some $m$. Iterating we see that $(xt^j)^{p^v}=x^{p^vm}t^{p^vj}$ for some $m$. Since $G_p$ has exponent $p^s$, $T/T^{p^{r+t}}$ has exponent $p^{r+t}$, and $r+t\geq s$, it follows that $(xt^j)^{p^{r+t}}=1$. This shows that the group has exponent ${p^{r+t}}$, so $H_p\to T/T^{p^{r+t}}$ is injective, and $\hat H$ is in fact $p$-torsion free.

This completes the proof of assembly for virtually cyclic groups.

\subhead 3.2 Torsion in rank one \endsubhead
The proof of Theorem 1.2.4 is given here because it is a simplified application of the methods used above. The theorem follows from:

\proclaim{3.2.1 Proposition}If $\Lambda$ is a subring of the rationals then
$$\Bbb H(E_{G,\Lambda-\t h.elem }/G;\Bbb K(p;R)\otimes\Lambda)\to  \Bbb H(E_{G,\t h.elem }/G;\Bbb K(p;R)\otimes \Lambda)$$
is an equivalence of spectra.\endproclaim
Here $\Bbb K \otimes \Lambda$ denotes the localized spectrum with homotopy groups $\pi_i(\Bbb K)\otimes \Lambda$. Localizing the coefficient system has the same effect as localizing the homology (i.e\. $\otimes \Lambda$ commutes with $\Bbb H$). 

The iterated homology identity (see \S5) expresses the map in 3.2.1 as induced by change of coefficients over $E_{G,\t h.elem }/G$, so it is sufficient to show the change of coefficients is an equivalence on stalks. We can also reduce to a single prime. This reduces to:

\proclaim{3.2.2 Lemma}  If\/ $G$ is a $q$-hyperelementary group and\/ $F$ is a free $G$--space then 
$$\Bbb H(E_{G,\t cyclic }/G;\Bbb K(p;R)\otimes\b Z [\frac1q])\to \Bbb  K(R[\pi_1 F/G])\otimes\b Z [\frac1q]$$
is an equivalence. \endproclaim
\subsubhead Proof of 3.2.2\endsubsubhead When $G$ is finite this is part of the classical theory (see Bartels-L\"uck \cite{BL}), so suppose $G$ is infinite. Suppose $F$ is a free $G$-space as usual. Let $A$ be a maximal normal infinite cyclic subgroup of $G$ and $Q$ the $p$-torsion quotient. In most cases we let $r$ be very large  and consider $G/A^{p^r}$. According to the classical theory  $K(R[\pi_1(F/G)])\otimes\b Z [\frac1p]$ is assembled from preimages in $G$ of {\it cyclic\/} subgroups $H\subset G/A^{p^r}$. In 3.1 we juggled two primes to find one so that $H_p$ is cyclic, and this is unnecessary here. The core of 3.1 similarly simplifies, and we go through it briefly.

Suppose $G$ is an infinite $p$-hyperelementary group with finite quotient $Q$ and let $G_p$ be a maximal $p$--torsion subgroup of $G$. Note $G_p\to Q$ is an injection. The proof divides into three cases:
\roster\item $Q$ acts nontrivially on $A$;
\item $G_p\to Q$ is an isomorphism and $Q$ acts trivially on $A$; and 
\item $G_p\to Q$ is not onto.
\endroster
As usual we proceed by induction on the size of $Q$, so the result is known for groups with smaller quotient. 
Case (1) is done last since it requires an extra wrinkle and logically depends on the other cases. 

Begin with case (2), so $G_p\to Q$ is an isomorphism. This gives a splitting of the sequence $1\to A\to G\to Q\to 1$, and $Q$ acts trivially on $A$ so $G=A\times Q$. The crystallographic quotient is the projection $A\times Q\to A$. Let $r$ be very large and consider assembly using cyclic subgroups $H\subset G/A^{p^r}$. The preimage $\hat H$ has finite quotient the image of $H\to Q$, so the result is known for $\hat H$ unless $H\to Q$ is onto. In particular if $Q$ is not cyclic this applies to all such $\Hat H$ and the theorem follows for $G$. Thus suppose $Q$ is cyclic, and denote its order by $p^s$. 

$G/A^{p^r}=A/A^{p^r}\times Q$, and $\{1\}\times Q$ is the (injective) image of the torsion in $G$. If $H$ intersects this trivially then $\hat H$ is torsion-free and therefore infinite cyclic. On the other hand it was shown in 3.1.3 that for $H$ to intersect $Q$ nontrivially it must lie in $A^{p^{r-s}}/A^{p^r}$. This means the image of $\hat H$ in $A$ lies in the image of the expansive map of index $p^{r-s}$. If $r$ is large enough this reduces $K$-theory to finite subgroups of $\hat H$, giving the theorem.

Next consider case (3), in which $G_p\to Q$ is not onto. As in 3.1.4 we show that if $r$ is large and $H\subset G/A^{p^r}$ is cyclic and maps onto $Q$ then $\hat H$ is torsion-free. 

The crystallographic quotient is cyclic, $G\to T$, so $G=G_p\rtimes_{\alpha} T$. If $p^r$ is greater than the order of $Q$ then $A^{p^r}=T^{p^{r+s}}$, where $p^s$ is the order of $Q/G_p$. In this case $H\subset G_p\rtimes T/T^{p^{r+s}}$. We only need consider $H$ that map onto $Q$, so map onto a nontrivial quotient of the cyclic $p$-group $T/T^{p^{r+s}}$, and consequently $H\to T/T^{p^{r+s}}$ must be onto. $\hat H$ intersects $G_p$ trivially if $H$ intersects $G_p\rtimes\{1\}$ trivially, and this is the kernel of the projection to $T/T^{p^{r+s}}$. Since $H$ is onto it intersects the kernel trivially if it has order $p^{r+s}$. But we saw in 3.1.4 that the semidirect product has exponent $p^{r+s}$ so this is a bound on the order of any cyclic subgroup. This concludes the proof in this case.

Finally consider case (1), with $p=2$ and $Q$ acting nontrivially on $A$. We assume $Q$ is cyclic, otherwise it can be reduced by cyclic assembly in $Q$ itself. Choose an odd prime $q$ and $r$ large, and consider $G/A^{q^r}$. This is 2-hyperelementary so $K\otimes  \b Z [\frac12]$ is assembled from cyclic subgroups. 

Let $G_2$ denote a maximal 2-torsion subgroup of $G$, then since the action is nontrivial the composition $G_2\to Q\to T_2=\text{Aut}(A)$ is onto. $Q$ is cyclic so $G_2\to Q$ must also be onto. This identifies $G$ as $A\rtimes G_2$, and $G/A^{q^r}$ as $A/A^{q^r}\rtimes G_2$. The crystallographic quotient is the dihedral group $A\rtimes T_2$. The image of a cyclic subgroup in $G/A^{q^r}$ is a cyclic subgroup of $A/A^{q^r}\rtimes T_2$. But the only cyclic subgroups of this that map onto the $T_2$ quotient are the order-two subgroups, and these have index $q^r$ in the dihedral quotient. This means the cyclic $H$ that map onto $Q$ have $\hat H$ with image in $A\rtimes T_2$ conjugate into $A^{q^r}\rtimes T_2$. But this is
 the image of the expanding map of index $q^r$. Consequently if $r$ is large enough then pulling back to $H$ reduces $K\otimes \b Z [\frac12]$ to finite cyclic subgroups. This completes the proof.

\subhead 3.3 Rank two \endsubhead
Here we treat the the rank--2 exceptions in Proposition 4.2.1. 
\proclaim{3.3.1 Proposition} Rank--two groups with an normal infinite cyclic subgroup satisfy hyperelementary assembly.\endproclaim
There is an extension $1\to A\to G\to Q\to 1$ with $A$ free abelian of rank 2. We write $A$ additively here to avoid confusion. Since there is an normal cyclic subgroup there is a summand of $A$ invariant invariant under the action of $Q$. Write $A=Z\oplus Z$ with the first summand invariant, then elements of $Q$ act by matrices of the form $\left(\smallmatrix  \pm1&a\\ 0& \pm1\endsmallmatrix\right)$. Changing the generator of the second summand changes the off-diagonal entry by an even number so we can standardize one such element to have off-diagonal entry 0 or 1. Then checking combinations that generate finite subgroups of $\t aut (Z\oplus Z)$ shows that $Q$ may have a $T_2$ summand generated by $\pm\left(\smallmatrix  1&1\\ 0& -1\endsmallmatrix\right)$ or $\pm\left(\smallmatrix  1&0\\ 0& -1\endsmallmatrix\right)$, and may have a $T_2$ summand generated by $\left(\smallmatrix  -1&0\\ 0& -1\endsmallmatrix\right)$. It is easy to determine the possible extensions and get a list of possibilities for $G$, but this is not needed here. 

As usual we consider hyperelementary subgroups of  quotients $G/A^n$. Let $H$ be such a subgroup and denote the preimage in $G$ by $\hat H$. The goal is to show that either $\hat H$ is known to satisfy hyperelementary assembly or there is a map to a rank-1 crystallographic group $G\to C$ that takes $\hat H$ into the image of an expansive map on $C$. In this case if estimates work out right then $K$-theory objects get pulled into the image of the finite-isotropy assembly on $C$ and therefore into the virtually-cyclic isotropy assembly on $G$. Previous results and finite hyperelementary induction then give the desired conclusion for $G$. 

The key point is the use of maps to rank-1 groups, and the proof divides into two cases depending on how many of these there are. If $A$ is inhomogeneous as a $Q$-module, i.e\. if an element acts by a matrix of the form $\pm\left(\smallmatrix  1&a\\ 0& -1\endsmallmatrix\right)$, then there are only two maps to rank-1 groups and not much can happen. The $Q$--homogeneous case is more delicate: there are infinitely many maps to rank-1 groups and most of them do not give the needed control. We must see that for each possible subgroup $H$ there is a rank-1 projection that works.

The proof in the homogeneous case comes from Farrell-Hsiang  \cite{FH2} and is based on the following lemma (c.f\.  \cite{FH2, Lemma 4.3}).

\proclaim{3.3.2 Lemma} Suppose $p$ is prime and\/ $C\subset (\b Z /p)^2$ is a nontrivial cyclic subgroup. Then there is a homomorphism $h\:\b Z ^2\to \b Z $ with $C$ the kernel of the mod\/ $p$ reduction of\/ $h$ and with norm $|h|<\sqrt{2p}$.  \endproclaim
\subsubhead Proof of the Lemma\endsubsubhead
If $C$ is contained in one of the standard factors then projection to the other factor works. If $C$ is not contained in a factor then it is generated by the equivalence class $[1,c]$ for some $c$ relatively prime to $p$. Since they are relatively prime there are $r,s$ with 
$rc +sp=1$. Now consider the integers $\{ra\mid 0\leq a\leq \sqrt{p}\}$. There are $\sqrt{p}$ of these with mod-$p$ representatives in the interval $[0,p-1]$ so there must be two  closer  together than $\sqrt{p}$. Denote the difference between  these two by $\alpha$ then there is $b$ so $|r\alpha+bp|<\sqrt{p}$. It is also the case that $|\alpha|<\sqrt{p}$ since it is the difference of nonnegative numbers less than $\sqrt{p}$. Now consider the homomorphism defined by $h(x,y)=\alpha x -(r\alpha+bp) y$. The estimates imply that $h$ has norm less than $\sqrt{2p}$. Finally we check that $C$ is the mod-$p$ kernel of $h$: 
$$h(1,c)=\alpha-r\alpha c-pbc=\alpha(1-rc)-pbc =(s-bc)p.$$
This completes the proof of the Lemma.

We can now prove the homogeneous cases of 3.3.1.
\proclaim{3.3.3 Lemma} $G=T\times T$ and $(T\times T)\rtimes T_2$ satisfy hyperelementary assembly.\endproclaim
\subsubhead Proof\endsubsubhead
Begin with $T\times T$. Suppose $p, q$ are large primes and $H\subset (T\times T)/(T\times T)^{pq}$ is $q$-hyperelementary. The $p$-torsion subgroup of $H$ is at most cyclic. According to Lemma 3.3.2 there is a homomorphism $h\:T\times T\to T$ so that the preimage $\hat H\subset G$ is taken to the subgroup $T^p$ and therefore factors through the $p$-expansive map $T\to T$.  Further, the norm of $h$ is less than $\sqrt{2p}$. 

We use this to get control. The description of $G$ as a product gives objects over the control space $S^1\times S^1$, and we measure path lengths there. Use $h$ to construct a map $S^1\times S^1\to S^1$. Applying $h$ increases path length by at most the norm of $h$, so no more than $\sqrt{2p}$, while transfer by the expansive map reduces it by a factor of $p$. The outcome is that path lengths in the image $S^1$  are changed by $\sqrt{2p}/p=\sqrt{2/p}$. Therefore we can get any specified reduction in lengths by choosing $p, q$ sufficiently large.

Now consider the twisted case and suppose $p$ is a large prime. If $H\subset G/(T\times T)^{p}$ is not onto the quotient of order 2 then it lies in the image of $T\times T$ and reduces to the abelian case.  Suppose therefore that $H$ is hyperelementary and onto this quotient. Then it is 2-hyperelementary and the $p$-torsion is at most cyclic. According to the lemma there is $h\:T\times T\to T$ with norm less than $\sqrt{2p}$ and with the $p$-torsion of $H$ in the kernel of the mod $p$ reduction. This extends to $\hat h\:(T\times T)\rtimes T_2\to T\rtimes T_2=D$, where $D$ is the dihedral group. The image in $D$ has index $p$ (and this is congruent to 1 mod 2)  so factors through the $p$-expansive homomorphism on $D$. 

Again we use this to get control. The description of $G$ provides a map to the quotient of $S^1\times S^1$ by an involution that fixes 4 points. The quotient is $S^2$ stratified with 4 ``singular'' points, and we use this to measure path lengths. The homomorphism $\hat h$ provides a stratified map $S^2\to I$. $\hat h$ increases path length by at most $\sqrt{2p}$ while transfer through the expansive map decreases it by $p$. Again there is a reliable reduction in length by $\sqrt{2/p}$, and we get control. 

This completes the proof of Lemma 3.3.2.
\subsubhead Proof of Proposition 3.3.1\endsubsubhead
Suppose $G$ is a group as specified in the Proposition. If it is not one of the groups of 3.3.3 then the action of $Q$ on $A$ has exactly two maximal invariant cyclic subgroups, and dividing by these gives exactly two homomorphisms of $G$ onto rank--1 crystallographic groups.  Furthermore the normalizations described above give crystallographic models in which the maps on quotients have derivatives of norm less than 2, so the maps increase path length by less than a factor of 2. 

Let $p$ be a large prime, as specified below, and consider $G/A^p$. Suppose $H$ is a hyperelementary subgroup, with preimage $\hat H\subset G$. If $\hat H$ has crystallographic quotient one of the groups in 3.3.3 then it satisfies hyperelementary assembly. If not then the action of the image  of $H\to Q$ on $\hat H\cap A$ again has two maximal invariant cyclic subgroups. These must be the same as the $Q$--invariant subgroups. $H$ must be 2-hyperelementary so the $p$--torsion subgroup $H\cap (A/A^p)$ is trivial or cyclic. It is invariant under the action of the image of $H\to Q$ so it lies in the the image of an invariant cyclic subgroup of $A$.  Let $G\to C$ be the map to a rank-1 group with this invariant subgroup in it's kernel. This implies that the image of $\hat H$ in $C$ has index $p$ and therefore lies in the image of the  $p$-expansive map. 

We now get control. The map to the crystallographic quotient of $C$ ($S^1$ or $I$) increases path length by a factor of at most 2, while pulling back through the expansive map decreases length by factor of $p$. This gives predictable improvement independent of choice of subgroup $H$ or rank-1 projection. Therefore choosing $p$ very large  pulls a large part of the $K$-theory into the controlled subspace and thus back through the assembly map. The assembly here is the finite-isotropy assembly on the rank-1 quotient, which by fibered reduction (\S2.3) lifts to  virtually-cyclic assembly in $G$. 

As explained at the beginning of the section this implies $G$ satisfies hyperelementary assembly, and the proof is complete. 

\head 4. Group theory \endhead
This section collects the group-theoretic results needed. Short direct proofs are given, avoiding for instance group cohomology, and largely drawn from work of Farkas \cite{F1, F2}. Several of the statements are minor but useful generalizations of the well-known versions.  

\subhead 4.1 The splitting group\endsubhead
There are standard maps of virtually abelian groups to semi-direct products. This treatment is adapted from Farcus \cite{F2}.
Suppose 
$$\xymatrix{1\ar[r]& A\ar[r] &G\ar[r]^{\pi} &Q\ar[r] &1}$$
is exact, with $A$ abelian and $Q$ finite. Denote the conjugation action of $Q$ on $A$ by $\alpha$. Choose  preimages $x_j\in \pi^{-1}(j)$ for each $j\in Q$, with $x_1=1$, and define a function $\beta\:G\to A$ by
$$\beta(y)=\prod_{j\in Q}y\,x_jx_{\pi(y)j}^{-1}.$$
Note that this depends on the choices $x_j$. 
\proclaim{4.1.1 Lemma} With the notation above $(\beta,\pi)\:G\to A\rtimes_{\alpha}Q$ is a homomorphism, and for $y\in A$, $\beta(y)=y^{|Q|}$.\endproclaim
The conjugation action of $Q$ on $A$ is independent of the preimage element chosen for the conjugation. Thus $\alpha(j)$ acts by $x_jax_j^{-1}$. The homomorphism statement is therefore equivalent to 
$$\beta(yz)=\beta(y)\,x_{\pi(y)}\beta(z)x_{\pi(y)}^{-1}.$$
To see this calculate:
$$\align \beta(yz)&=\prod_{j}yzx_jx_{\pi(yz)j}^{-1}\\
&=\prod_jy(x_jx_{\pi(y)j}^{-1}x_{\pi(y)j}x_j^{-1})(x_{\pi(y)}^{-1}x_{\pi(y)})zx_j(x_{\pi(z)j}^{-1}x_{\pi(y)}^{-1}x_{\pi(y)}x_{\pi(z)j})x_{\pi(yz)j}^{-1}\endalign$$
Parse this into terms in $A$ and use the fact that $A$ is abelian to rearrange to
$$\prod_jyx_jx_{\pi(y)j}^{-1}\prod_j x_{\pi(y)j}x_{j}^{-1}x_{\pi(y)}^{-1}
\prod_j x_{\pi(y)}zx_{j}x_{\pi(y)j}^{-1}x_{\pi(y)}^{-1}\prod_jx_{\pi(y)}x_{\pi(z)j}x_{\pi(y)\pi(z)j}^{-1}$$
The first term is $\beta(y)$ and the third is $x_{\pi(y)}\beta(z)x_{\pi(y)}^{-1}$. Reindex the fourth term with $k=\pi(z)j$, then this is seen to be the inverse of the second term so these cancel.

Finally if $y\in A$, so $\pi(y)=1$, then $\beta(y)=\prod_{j\in Q}yx_jx_j^{-1}=y^{|Q|}$.

\subhead 4.2 Crystallographic quotients\endsubhead
A group is crystallographic if it is an extension of a free abelian group by a finite group, and the conjugation action gives an injection of the finite group into the automorphisms of the abelian subgroup, see \cite{F2}. 

\proclaim{4.2.1 Lemma} A finitely generated virtually abelian group maps onto  a crystallographic group of  the same rank.\endproclaim
Denote the group by 
$$\xymatrix{1\ar[r]& A\ar[r] &G\ar[r]^{\pi} &Q\ar[r] &1}$$
as above, then the finite generation hypothesis implies $A$ is finitely generated. Denote by $\bar{\alpha}$ the action of $Q$ on the torsion-free quotient:
$$Q\to \t Aut (A)\to \t Aut (A/\t torsion ).$$
Then the quotient maps define a morphism of the splitting group of 4.1,
$$A\rtimes_{\alpha}Q\to (A/\t torsion )\rtimes_{\bar{\alpha}}(\t image \,\bar{\alpha}).$$
The image of $G$ in the second group is crystallographic, and it has finite index so has the same rank as $G$.

\subhead 4.3 Expansive maps \endsubhead
These provide the key geometric input for the main theorem. Again the treatment is adapted from Farcus \cite{F2}. 

For $r\in Z$ define an endomorphism of the splitting group $A\rtimes Q$ by $(\wedge r,\t id )(a,j)=(a^r,j)$. 
 Let $\beta$ be the function defined in 4.1, and recall that it depends on choices.
\proclaim{4.3.1 Lemma} For any $n\in Z$ the diagram commutes:
$$\xymatrix{G\ar[rr]^{\beta^n\cdot \t id }\ar[d]^{(\beta,\pi)}&&G\ar[d]^{(\beta,\pi)}\\
A\rtimes Q\ar[rr]^{(\wedge(1+n\,|Q|),\t id )}&&A\rtimes Q}$$\endproclaim
The top map is $g\mapsto \beta(g)^ng$, and  can be checked to be a homomorphism.  Commutativity follows easily from the fact that $\beta(y)=y^{|Q|}$ for $y\in A$. The lemma motivates calling the top homomomorphism an {\it expansive map with expansion $1+n\,|Q|$}.

\proclaim{4.3.2 Corollary} Suppose $H\subset G$ is a subgroup with no $r$-torsion, and $r\equiv 1$ mod $|Q|$. Then the inclusion of $H$ factors through the expansive map of expansion $r$ if and only if the image of $(H\,A^r)/A^r$  in the quotient of the splitting group $(A/A^r)\rtimes Q$ lies in $1\rtimes Q$.\endproclaim
This is because the image of the map $(\wedge(1+n\,|Q|, \t id )$ on the splitting group is exactly the preimage of $1\rtimes Q$ in the quotient.

\proclaim{4.3.3 Corollary} Suppose $H\subset G$ is a subgroup with no $r$-torsion, and $r\equiv 1$ mod $|Q|$. Then the following are equivalent:
\roster \item $H$ is conjugate to a subgroup that factors  through the expansive map of expansion $r$;
\item $(H\,A^r)/A^r \subset (A/A^r)\rtimes Q$ is $|Q|$-torsion;
\item $(H\,A^r)/A^r\to Q$ is an injection.
\endroster
\endproclaim
This follows easily from 4.3.2.

\head 5. The Farrell-Jones component\endhead
This section provides proofs of statements in \S1.4 about the Farrell-Jones component of virtually-cyclic isotropy homology. The main ingredient is the elaborate functoriality of this version of $K$-theory, used to obtain interesting general results from special cases. This is illustrated in 5.1 with a version of the ``Novikov conjecture''.  In 5.2 it is used to give a homology description of the Farrell-Jones  part of hyperelementary or virtually-cyclic isotropy $K$ homology. This is shown to be a split summand in 5.3. 

\subhead 5.1 The ``Novikov conjecture''\endsubhead
As a warmup for the next section we show:
 \proclaim{5.1.1 Proposition} If $G$ satisfies hyperelementary assembly then the assembly map 
$$\Bbb H(BG;\Bbb K(R))\to \Bbb K(R[G])$$
is a rational split injection. \endproclaim
This begins with a trivial general construction. Suppose $p\:E\to X$ is a stratified system of fibrations. We can regard the horizontal maps in the diagram 
$$\xymatrix{E\ar[d]^{\text{id}}\ar[r]^{\text{id}}&E\ar[d]^p\ar[r]^p&X\ar[d]^{\text{id}}\\
E\ar[r]^p&X\ar[r]^{\text{id}}&X}$$
 as morphisms of the vertical maps thought of as stratified systems of fibrations. 
Applying $\Bbb K$ fiberwise gives  morphisms of spectral cosheaves:
$$\xymatrix{E\times{\Bbb K(R)}\ar[d]\ar[r]&{\Bbb K}(p;R)\ar[d]\ar[r]&X\times{\Bbb K(R)}\ar[d]\\
E\ar[r]^p&X\ar[r]^{\text{id}}&X}$$
(The first and last cosheaves  are constant.)
These induce morphisms on homology,
$$ \Bbb H(E;\Bbb K(R))\to \Bbb H(X;\Bbb K(p,R))\to \Bbb H(X;\Bbb K(R)).$$
In the particular case that $E\to X$ is  $E_G/G\to  E_{G,\text{finite}}/G$ this gives
$$ \Bbb H(E_G/G;\Bbb K(R))\to \Bbb H(E_{G,\text{finite}}/G;\Bbb K(p,R))\to \Bbb H(E_{G,\text{finite}}/G;\Bbb K(R)).\tag{5.1.2}$$
We  refine this using the iterated homology identity. This identity expresses the homology of the domain of a map as homology of the range with coefficients in homology of the point inverses, and is a spectrum version of the Leray spectral sequence. Returning to the abstract $p\:E\to X$ we write 
$$\Bbb H(E;\Bbb K(R))=\Bbb H(X;\Bbb H(p;\Bbb K(R))).$$
$\Bbb H(p;\Bbb K(R))$ is the system over $X$ with stalk over $x$ the spectrum $\Bbb H(p^{-1}(x);\Bbb K(R))$. Now we can express the morphism coming from the diagram of maps as induced by morphisms of coefficient systems over $X$. 

In the group action application, the maps in (5.1.2) are induced by morphisms of systems over $E_{G,\text{finite}}/G$. Over a point with isotropy group $J$ the map of stalks is
$$\Bbb H(E_J/J;\Bbb K(R))\to \Bbb K(R[J])\to \Bbb K(R).$$
The composition is the map induced on $\Bbb K$ homology by the map $E_J/J\to \t pt $. But $J$ is a finite group, so this is a rational equivalence (in fact a $1/|J|$ equivalence) and induces a rational equivalence on homology. 

Returning to diagram 5.1.2, we have seen that the composition of the two maps is induced by a rational equivalence of coefficient systems. The induced map on homology is therefore a rational equivalence, and in particular the first map is a rational split injection.

We will see below that the finite-isotropy homology injects into the virtually-cyclic isotropy homology, so if this is the same as $K$-theory the proposition follows. This completes the proof of 5.1.1.

The argument reveals that it is only necessary to invert primes dividing orders of finite subgroups of $G$, but it is almost certainly false if these are not inverted. It also shows that the image of the constant-coefficient assembly is disjoint from the Farrell-Jones component of $K$-theory. One conclusion is that the Novikov conjecture is not a particularly good way to approach $K$-theory of groups with torsion. 

\subsubhead 5.1.3 Basepoints\endsubsubhead
A caution about basepoints and naturality may be appropriate here. The stalks in the coefficient system ${\Bbb K}(p;R)\to E_{G,\text{finite}}/G$ used in the middle in 5.1.2 can be identified as ${\Bbb K}(R[Q])$ for finite subgroups $Q\subset G$. The second map in 1.4.1 is obtained by projecting to the trivial group $Q\to 1$ in each fiber, then applying $\Bbb K$ fiberwise. It is tempting to say ``apply $\Bbb K$ fiberwise to the {\it injection\/} $1\to Q$ in each fiber'', and therefore get a map of coefficient systems $(E_{G,\text{finite}}/G)\times{\Bbb K(R)}\to{\Bbb K}(p;R)$ over $E_{G,\text{finite}}/G$. This would give an {\it integral\/} splitting of the second map in 5.1.2. The problem with this argument is that identifying a stalk of the coefficient system with  ${\Bbb K}(R[Q])$ requires choice of a basepoint. The $K$-theory map induced by $Q\to 1$ comes geometrically from a map to a point and does not depend on the basepoint. However the inclusion $1\to Q$ corresponds to including the basepoint. The map $E_G/G\to  E_{G,\text{finite}}/G$ does not have a section so there is no continuous choice of basepoints for fibers and this argument does not work.

\subhead 5.2 Proof of Proposition 1.3.2\endsubhead 
The objective is to express the Farrell-Jones component as homology with coefficents in $\Bbb F\Bbb J$ applied fiberwise. 

$\Bbb F\Bbb J(G;F,R)$ is defined as the cofiber of the map of spectra
$$\Bbb H(E_{G,\text{finite}}/G;\Bbb K(p_{\text{finite}},R))\to \Bbb H(E_{G,\text{v.cyclic}}/G;\Bbb K(p_{\text{v.cyclic}},R)).$$
Rewrite the first term as a homology of $E_{G,\text{v.cyclic}}/G$ using the iterated homology identity, then this map is induced by a morphism of coefficient systems 
$$\Bbb H(q;\Bbb K(p_{\text{finite}},R))\to \Bbb K(p_{\text{v.cyclic}},R)$$
where $q$ is the projection $(E_{G,\text{finite}}\times E_{G,\text{v.cyclic}})/G\to E_{G,\text{v.cyclic}}/G$. Over a point with a virtually cyclic isotropy group $J$ the map of stalks is equivalent to 
$$\Bbb H(E_{J,\text{finite}}/J;\Bbb K(p_{\text{finite}},R))\to \Bbb K(R[\pi_1(F/J)]).$$
Since $J$ is virtually cyclic $E_{J,\text{v.cyc}}/J$ is contractible, so we can think of the target of this map as homology of a point, $\Bbb H(E_{J,\text{v.cyc}}/J;\Bbb K(p_{\text{v.cyc}},R))$. Using this we recognize the cofiber of this map of stalks as (by definition) $\Bbb F\Bbb J(J;F,R)$. In other words the cofiber of the morphism of coefficient systems is the coefficient system obtained by applying $\Bbb F\Bbb J$ fiberwise. 

The final ingredient is the observation that the cofiber of an induced map on homology is the homology with coefficients in the stalkwise cofiber of the morphism of coefficient systems. 
\subhead 5.3 Relaxation\endsubhead 
Here we prove Proposition 1.3.3 asserting that the Farrell-Jones component is a split summand of the virtually-cyclic isotropy homology.
\proclaim{5.3.1 Theorem} For any $G$ and free $G$-space $F$ there is a natural left inverse for the map
$$\Bbb H(E_{G,\text{finite}}/G;\Bbb K(p_{\text{finite}},R))\to \Bbb H(E_{G,\text{v.cyc}}/G;\Bbb K(p_{\text{v.cyc}},R)).$$
\endproclaim
This extends a result of Bartels \cite{B}. We refer to the left inverse as  ``relaxation'' because it is an extension of an operation with that name developed in the context of Whitehead groups of twisted polynomial extensions and fibering manifolds over  a circle, by Farrell, Siebenmann and others.

Recall that in 5.2 we identified the map in this theorem as being induced on homology by a morphism of coefficient systems over $E_{G,\text{v.cyc}}/G$. To split the homology map it is sufficient to split the morphism of coefficient systems. Further over each point  the map of stalks is again a map of the form of 5.3.2, but applied to a virtually cyclic subgroup of $G$. Therefore to prove the theorem in general it is sufficient to prove it for virtually cyclic groups, provided the inverse map is natural  enough to allow it to  be applied ``continuously'' fiberwise  over a base space. Technically the extension uses simplicial functors \cite{Q1 \S6} rather than topological continuity. We state the group version because it is easier to understand, even though basepoint problems keep it from being sufficiently natural, c.f.\ 5.1.3.  A more natural formulation is given at the beginning of the proof.

\proclaim {5.3.2 Lemma} Suppose $G$ is infinite virtually cyclic, $F$ is a free $G$-space, and $G\to C$ is the crystallographic quotient. Then the assembly map
$$\Bbb H(\b R /C;\Bbb K(p,R))\to \Bbb K(R[\pi_1(F/G)])$$
has a natural left inverse up to homotopy.\endproclaim
  Recall $C$ is either infinite cyclic or dihedral and $\b R /C$ is $S^1$ or $I$ respectively. The induced action of $G$ on $\b R $ is a model for $E_{G,\text{finite}}$. The coefficients in the homology group come (as usual) from $p\:(F\times \b R )/G\to \b R /G$.
  
For the more natural formulation we use the notation $\Bbb K(F/G;R)$ as shorthand for the space of $K$ data in $F/G$ without control (or equivalently, control over a point). According to Theorem 2.2.1 of \cite{Q1} a choice of basepoint determines an equivalence of this with the Quillen $K$-spectrum $\Bbb K(R[\pi_1(F/G)])$ on the right in the Lemma. The statement we prove is: there is a canonical spectrum $\Bbb J$ and a commutative diagram  
$$\xymatrix{{\Bbb H}(\b R /C;\Bbb K(p,R))\ar[r]^-A\ar[dr] &{\Bbb K}(F/G;R)\ar[d]\\
&{\Bbb J}}\tag{5.3.3}$$
and the diagonal map is an equivalence.  The canonical homotopy class of the Lemma is obtained by composing the  map $\Bbb K\to \Bbb J$ with a homotopy inverse for the diagonal. The result is ``homotopy-everything'' well defined because homotopy inverses are, but we can avoid spelling this out by including the intermediate space $\Bbb J$ in the statement. Applying this fiberwise in a stratified situation gives spectral cosheaves with fibers $\Bbb H$, $\Bbb K$ and $\Bbb J$, and morphisms as in 5.3.3 between them. Applying homology to a fiberwise homotopy equivalence of coefficient systems gives a homotopy equivalence on homology, and the inverse map needed for the global version 5.3.1 is obtained from a homotopy inverse for this equivalence. 

We begin the proof of Lemma 5.3.2 with $G$ with cyclic quotient: $C\simeq T$. $G$ is a semidirect product $P\rtimes_{\alpha}T$ where $P$ is the torsion subgroup and $T$ acts by an  automorphism $\alpha$. The map $p$ and the corresponding coefficient system are bundles (the quotient $S^1$ is all one stratum). More explicitly, $\alpha$ induces an automorphism of the classifying space $E_P\to E_P$ and $p$ is the projection of the corresponding $(E_P\times F)/P\simeq F/P$ bundle over $S^1$. The coefficient system is the $\Bbb K(R[\pi_1(F/P)])$ bundle over $S^1$ now described either by applying $\Bbb K$ fiberwise to the $(E_P\times F)/P$ bundle or by using the automorphism  induced by $\alpha$. Let $\hat p$ denote the pullback  over the universal cover $\b R \to S^1$. 

The comparison with $K$-theory is made using exactness in locally finite homology. We use the universal cover for this rather than open sets, both to make it more natural and to fit with the $K$-theory construction. Pulling back to the universal cover maps ordinary homology to locally finite homology:
$$\Bbb H(S^1; \Bbb K(p,R))\to \Bbb H^{lf}(\b R ; \Bbb K(\hat p,R)).$$
The deck transformation $+1\:\b R \to \b R $ is covered by a morphism $\hat \alpha$ of the coefficient system. If we choose a basepoint in $\b R $ then we get a trivialization of the coefficient system, $\Bbb K(\hat p;R)\simeq \b R \times \Bbb K(R[P])$. In these coordinates the deck-transformation morphism is $+1$ in the $\b R $ coordinate and $\Bbb K(\alpha)$ in the $\Bbb K(R[P])$ coordinate. On the pullback from $S^1$ there is a canonical homotopy between the identity and the map induced by this morphism. Since we are working with spectra we can subtract and get a canonical nullhomotopy of the difference. A straightforward deduction from the exact sequence of a pair is that this nullhomotopy makes the sequence
$$\xymatrix{{\Bbb H}(S^1; \Bbb K(p,R))\ar[r]& {\Bbb H}^{lf}(\b R ; \Bbb K(\hat p,R))\ar[r]^{1-\alpha^*}&{\Bbb H}^{lf}(\b R ; \Bbb K(\hat p,R))}\tag{5.3.4}$$
a homotopy fibration. More precisely let $\Bbb I$ denote the homotopy fiber of the second map, then the nullhomotopy induces a map $ {\Bbb H}(S^1; \Bbb K(p,R))\to \Bbb I$ and this is an equivalence of spectra. 

Now do a similar construction in $K$-theory. Begin with $\Bbb K(F/G;R)$. Pullback to the covering space induced by the map to $S^1$ goes to a space we denote temporarily by  $\Bbb K^{lf}(\b R ;F/P;R,\text{bdd})$.  These are $K$-type objects with data in $\b R \times F/P$ that may be infinite but have locally finite image in $\b R $ and lengths of paths and homotopies measured in $\b R $ are bounded. The proof below works without the bounded hypothesis but further argument would be needed and bounded versions of $K$-theory are more familiar, c.f\. \cite{P}. 

The endomorphism of  $\Bbb K^{lf}(\b R ;F/P;R,\text{bdd})$ induced by the deck transformation in the cover is canonically homotopic to the identity on the image of $\Bbb K(F/G;R)$. It also is induced by the endomorphism $\alpha$ of $F/P$. This gives a nullhomotopy of the composition of the maps in the diagram
$$\xymatrix{{\Bbb K}(F/G;R)\ar[r]& {\Bbb K}^{lf}(\b R ;F/P;R,\text{bdd})\ar[r]^{1-\alpha^*}&{\Bbb K}^{lf}(\b R ;F/P;R,\text{bdd})}\tag{5.3.5}$$
The $\Bbb J$ we seek is the homotopy fiber of the second map in this diagram. The nullhomotopy defines a map  
${\Bbb K}(F/G;R)\to \Bbb J$.

There are assembly maps from diagram 5.3.4 to 5.3.5:
 $$\xymatrix{{\Bbb H}(S^1; \Bbb K(p,R))\ar[r]\ar[d]^{A}& {\Bbb H}^{lf}(\b R ; \Bbb K(\hat p,R))\ar[r]^{1-\alpha^*}\ar[d]^{A}&{\Bbb H}^{lf}(\b R ; \Bbb K(\hat p,R))\ar[d]^{A}\\
{ \Bbb K}(F/G;R)\ar[r]& {\Bbb K}^{lf}(\b R ;F/P;R,\text{bdd})\ar[r]^{1-\alpha^*}&{\Bbb K}^{lf}(\b R ;F/P;R,\text{bdd}).}\tag{5.3.6}$$
Naturality of assembly maps implies that this diagram commutes, and induces a map of homotopy fibers of the right-hand maps $\Bbb I\to \Bbb J$. Moreover the right two assembly maps in this diagram are equivalences because over $\b R $ any number is a threshold for stability. To see this last, we know there is some threshold $\epsilon$. Given any $r>0$ we can compress $\b R $ to shrink things of size $<r$ to size $<\epsilon$. Apply stability and then decompress the result. The conclusion from this is that the induced map of fibers  $\Bbb I\to \Bbb J$ is also an equivalence. 

The nullhomotopies of the compositions in 5.3.4 and 5.3.5 are compatible in 5.3.6 because they have essentially the same definition. Compatibility of the homotopies means the diagram of maps to fibers commutes:
$$\xymatrix{{\Bbb H}(S^1; \Bbb K(p,R))\ar[r]^A\ar[d]^{\simeq}& { \Bbb K}(F/G;R)\ar[d]\\
{ \Bbb I}\ar[r]^{\simeq}&{ \Bbb J}}.$$
This is the conclusion we need for the applications of the lemma. 

We finish the proof of Lemma 5.3.2 with a similar construction for $G$ with dihedral image  $G\to D$. Think of $D$ as acting on $\b R $ so the generator of the index-two subgroup acts by translation by 2. 
The quotient $\b R /D$ is then the unit interval. The control map $p\:(F\times \b R )/G\to I$ induces a stratification with two strata: the endpoints and the interior. The fundamental group over the interior is index two in the groups over the endpoints, and $\pi_1(F/G)$ is the free product of the endpoint groups amalgamated over the interior group. We use a Meyer-Vietoris sequence in locally finite homology to analyze this. Write the interval $[0,1]$ as $[0,1)\cup (0,1]$ then  restrictions $[0,1]\to [0,1)\cup (0,1]\to (0,1)$ give a long exact sequence of groups, or more precisely a homotopy fibration of spectra. As in the cyclic case we rewrite this using covers instead of subsets. Let $T_2^{(0)}$ be the order-two subgroup of $D$ whose nontrivial element is the involution fixing 0, and let $G_2^{(0)}$ denote the inverse image of this in $G$. Then there is a diagram 
$$\xymatrix{(F\times \b R )/G^{(0)}_2\ar[r]^{q_0}\ar[d]^{p^{(0)}}&  (F\times \b R )/G\ar[d]^p\\
{ \b R }/T^{(0)}_2\ar[r] & {\b R }/D.}\tag{5.3.7}$$
The bottom map  is the
 folding map  $[0,\infty)\to [0,1]$ and is locally a covering map except at positive integers where it is a 2-fold branched cover branched over an endpoint. 
The map $q_0$ across the top is a covering map so the diagram induces a pullback map in locally-finite homology,
$$q_0^*\: \Bbb H([0,1];\Bbb K(p;R))\to  \Bbb H([0,\infty);\Bbb K(p^{(0)};R)).$$

The left vertical map in 5.3.7 has locally constant fundamental group over $[0,\infty)$ except at the remaining basepoint. This means $[0,1)\to [0,\infty)$ induces an equivalence on locally-finite homology with $\Bbb K(p;R)$ coefficients, and shows the pullback $[0,1]\to [0,\infty)$ is the same up to homotopy as the restriction.

 There is a similar construction at the other end of $[0,1]$. Let $T_2^{(1)}\subset D$ be generated by the involution fixing 1, and let $G_2^{(1)}$ denote the inverse image of this in $G$. $\b R /T_2^{(1)}$ is $(-\infty,1]$ and the map to $[0,1]$ is again a folding map and is covered by the covering map $q_1\:(F\times \b R )/G^{(1)}_2\to (F\times \b R )/G$. The induced pullback in locally finite homology is equivalent to the restriction induced by $(0,1]\to [0,1]$. 
 
 Both of these covering constructions have a further 2-fold cover branched at the remaining endpoint, ending up over $(-\infty,\infty)$. The compositions $[0,1]\to  (-\infty,\infty)$ is again a pullback and is independent of which half-infinite interval it goes through. The difference between the compositions is therefore trivial, so there is a trivialization of the composition in the following sequence:
$$\xymatrix{{\Bbb H}(I;\Bbb K)\ar[r]^-{( q_0^*, q_1^*)}&{\Bbb H}^{lf}([0,\infty);\Bbb K)\times{\Bbb H}^{lf}((-\infty,1];\Bbb K)\ar[r]^-{q_1^*-q_2^*}&{\Bbb H}^{lf}((-\infty,\infty);\Bbb K). }\tag{5.3.8}$$
Denote the homotopy fiber of the second map by $\Bbb I$, then identification of this with the Meyer-Vietoris sequence shows the induced map   ${\Bbb H}(I;\Bbb K))\to \Bbb I$ is a homotopy equivalence. 

Now we do a similar construction in $K$-theory. Pullback to the covering takes $K$ objects in $(F\times \b R )/G\to [0,1]$ to $K$ objects in $(F\times \b R )/G_2^{(0)}\to [0,\infty)$ that are locally finite and bounded in the $[0,\infty)$ coordinate. This gives maps
$$\xymatrix{{\Bbb K}(F/G;R)\ar[r]^-{( q_0^*, q_1^*)}&{\Bbb K}^{lf}([0,\infty);p^{(0)},\t bdd )\times{\Bbb K}^{lf}((-\infty,1];p^{(1)},\t bdd  )\\
\ar[r]^-{q_1^*-q_2^*}&{\Bbb K}^{lf}((-\infty,\infty);\hat p,\t bdd ). }\tag{5.3.9}$$
  Again the composition is trivial because pullback to $(-\infty,\infty)$ is independent of which half-infinite pullback is used as an intermediate step. Denote the homotopy fiber of the second map by $\Bbb J$, then the triviality induces a map ${\Bbb K}(F/G;R)\to \Bbb J$. 
  
There are assembly maps taking the homology sequence 5.3.8 to 5.3.9, giving a commutative diagram by naturality. This induces a map of homotopy fibers $\Bbb I\to \Bbb J$, and since the trivializations of compositions also commute there is a commutative diagram 
$$\xymatrix{{ \Bbb H}([0,1];\Bbb K(p;R))\ar[r]^A\ar[d]&{\Bbb K}(F/G;R)\ar[d]\\
{ \Bbb I}\ar[r]&{ \Bbb J}.}$$
We know from the Meyer-Vietoris argument that the left vertical map is a homotopy equivalence. It only remains to see that the bottom map is a homotopy equivalence, and for this we need that the assembly maps connecting the center and right spaces in diagrams 5.3.8 and 5.3.9 are homotopy equivalences. The right one was already discussed in the cyclic-quotient case, and the center maps $ \Bbb H([0,\infty);\Bbb K(p^{(0)};R))\to{\Bbb K}^{lf}([0,\infty);p^{(0)},\t bdd )$ are equivalences for the same reason: any number is a stability threshold for the control map 
$$\xymatrix@1{(F\times \b R )/G^{(0)}_2\ar[r]^-{p^{(0)}}&  { \b R }/T^{(0)}_2\simeq[0,\infty)}.$$
This is easily seen using compression toward the basepoint.

This completes the proof of Lemma 5.3.2, and therefore 5.3.1. 
\head 6. Universal spaces\endhead
Universal $G$-spaces with specified isotropy are used extensively in this theory. L\"uck \cite{L} has given a survey with general constructions and many examples. Here we give a simple model that makes some of the structure particularly transparent. 

\subhead 6.1 Definition\endsubhead Suppose $\Cal G$ is a collection of isomorphism classes of groups and $G$ is a group. Define $E_{G,\Cal G}$ to be the simplicial set with $n$-simplices sequences $(C_0,C_1,\dots,C_n)$ where
\roster\item Each $C_i$ is a left coset of a subgroup of $G$; 
\item for each $i$, $C_iC^{-1}_i \subset  C_{i+1}C^{-1}_{i+1}$; and
\item $C_iC^{-1}_i $ is in the class $\Cal G$.
\endroster
To clarify (2) and (3) note that if $C_i$ is a coset $gH$ then $C_iC^{-1}_i =gHg^{-1}$ is a subgroup of $G$. Face and degeneracy maps are defined by omitting or duplicating entries respectively. $G$ acts on this diagonally: $g(C_0,C_0,\dots,C_n)=(gC_0,gC_1,\dots,gC_n)$. The spread-out nature of the action means this does not give a nice model for the quotient, and this is probably why it is not more widely used.

The subgroup of $G$ fixing a simplex $(C_*)$ is $C_0C_0^{-1}$. Similarly the fixed complex of a subgroup $H\subset G$ is the set of simplices $(C_*)$ with $C_0C_0^{-1}=H$. 
\subhead 6.2 Universality\endsubhead The universal space is  supposed to receive a $G$-map, well-defined up to homotopy, from every $G$-space with isotropy in $\Cal G$. We briefly describe these maps.

Suppose $X$ is a discrete $G$-set. Choose a subset $\bar X$ containing one point in each orbit, then the action gives a surjection $G\times \bar X\to X$. Point inverses determine cosets of subgroups: if $y=gx$ for $x\in \bar X$ then the inverse image of $y$ is $gH\times\{x\}$, where $H$ is the isotropy group of $x$. The assignment $y\mapsto (gH)$ defines a $G$-map from $X$ to the 0-simplices of $E_{G,\Cal G}$ provided $\Cal G$ contains the isotropy groups of $X$. 

Now suppose $X$ is a simplicial $G$-set. Choices for the 0-simplices as above defines a map $X^{(0)}\to E_{G,\Cal G}$. This extends to a simplicial map by taking a simplex to the sequence of images of its vertices. Note that for this to work we must define ``simpicial $G$-set'' to mean if an element of $G$ takes a simplex to itself then it fixes the simplex, and if it fixes a vertex of a simplex then it fixes vertices later in the ordering. The second condition can be arranged by subdividing and reordering, if it is not initially satisfied. Different basepoint choices give homotopic maps.

This construction directly gives the maps required for the universal property. It is also easy to check the contractibility criterion. The subcomplex fixed by $H\subset G$ is the set of $(C_*)$ with $C_0C_0^{-1}=H$. We can extend the identity map of this to itself to a map of the cone: the cone on a simplex $(C_*)$ goes to $(H,C_0,\dots,C_n)$. This gives a contraction that is equivariant with respect to the action of the normalizer of $H$.
\subhead 6.3 Unions\endsubhead Suppose $\Cal G$, $\Cal H$ are collections of isomorphism classes of groups. We say $\Cal G$ is {\it closed relative to } $\Cal H$ if a group that is a subgroup of an element of $\Cal G$ and is in $\Cal H$ must in fact be in $\Cal G$. 
\proclaim{Lemma 6.3.1} Suppose $\Cal G$, $\Cal H$ are collections of isomorphism classes closed with respect to each other. Then for a group $G$
\roster\item $E_{G,\Cal H\cup\Cal G}=E_{G,\Cal H}\cup E_{G,\Cal G}$ and
\item $E_{G,\Cal H\cap\Cal G}=E_{G,\Cal H}\cap E_{G,\Cal G}$
\endroster\endproclaim
To see (1) consider a simplex $(C_0,\dots, C_n)$ in $E_{G,\Cal H\cup\Cal G}$. If $C_nC_n^{-1}$ is in $\Cal G$ then the closure property shows all $C_iC_i^{-1}$ are in $\Cal G$, so the simplex is in $E_{G,\Cal G}$. The alternative is that it is in $E_{G,\Cal H}$. Statement (2) is similar.

\subsubhead 6.3.2 Example\endsubsubhead If $p$, $q$ are different primes, then 
$$E_{G,p-\t h.elem }\cup E_{G,q-\t h.elem }=E_{G,p \text{ or }q-\t h.elem }$$
$$E_{G,p-\t h.elem }\cap E_{G,q-\t h.elem }=E_{G, \text{cyclic} }.$$

\subhead 6.4 Joins\endsubhead Suppose $X$ and $Y$ are simplicial sets. Formally adjoin a simplex $\emptyset$ of degree $-1$ to each of them. Define a face map on 0-simplices to take everything to $\emptyset$. The {\it join\/} is then defined to be the set with $n$-simplices 
$$(X*Y)^{(n)}= \cup_{i+j=n-1}X^{(i)}\times Y^{(j)}.$$
Face maps on $\delta^i\times \sigma^j$ are defined by $\partial_k\times \t id $ for $k\leq i$ and $\t id \times \partial_{k-i-1}$ for $k>i$. Degeneracies are defined similarly. 
Note that the formal $-1$--simplices give inclusions $X\simeq X\times\{\emptyset\}\subset X*Y$ and $\{\emptyset\}\times Y\subset X*Y$. 

This connects with universal spaces as follows:
\proclaim{Lemma 6.4.1} Suppose $\Cal H\subset \Cal G$ are collections of group isomorphism classes, $\Cal H$ is closed relative to $\Cal G$, and $G$ is a group. Then there is a natural $G$-equivariant inclusion
$$E_{G,\Cal G}\subset E_{G,\Cal H}*E_{G,\Cal G-\Cal H}.$$
\endproclaim
Suppose $(C_0,C_1,\dots,C_n)$ is a simplex of $E_{G,\Cal G}$. If there is $i$ with  $C_iC_i^{-1}$  in $\Cal H$ then we let $j$ be the maximal such $i$, and otherwise set $j=-1$. The closure property of $\Cal H$ implies that $C_iC_i^{-1}$ is in $\Cal H$ if $i\leq j$ and is in $\Cal G-\Cal H$ if $i>j$. Therefore $(C_0,\dots,C_j)$ is a simplex of $E_{G,\Cal H}$ (the virtual $-1$ simplex $(\emptyset)$ if $j=-1$), and similarly $(C_{j+1},\dots,C_n)$ is a simplex of $E_{G,\Cal G-\Cal H}$. The inclusion of the lemma takes $(C_*)$ to the simplex $(C_0,\dots,C_j)*(C_{j+1},\dots,C_n)$ in the join.  This also identifies the image as joins $(C_0,\dots C_k)*(D_0,\dots, D_n)$ with $C_kC_k^{-1}\subset D_0D_0^{-1}$. This can be elaborated to get a description of $E_{G,\Cal G}$ as a pushout, but here we use something weaker.

There is a natural homeomorphism of the realization  $|X*Y|$  to the topological join $X\times [0,1]\times Y/\sim$, where $(x,t,y)\sim (x',t',y')$ if either they are equal, or $t=t'=0$  and $x=x'$, or $t=t'=1$ and $y=y'$. The natural homeomorphism in general follows from the case where $X$ and $Y$ are simplices.
\proclaim{Lemma 6.4.2} Suppose $Z\subset X*Y$ is a subcomplex containing  $Y$. Then the complement of $|X|$ in the realization $|Z|$ deformation retracts to $|Y|$. \endproclaim
The complement of the $|X|$ in the topological join is the open mapping cylinder of the projection $|X|\times|Y|\to |Y|$, and the deformation is the obvious one using the mapping cylinder coordinate. The deformation preserves simplices not in $|X|$ so gives a deformation of any subcomplex.  

\subsubhead 6.4.3 Example\endsubsubhead 
The class of cyclic groups is closed relative to all groups, so in particular relative to hyperelementary groups. We therefore get an embedding of $E_{G,\t h.elem }$ as a subcomplex of the join $E_{G,\t cyclic }*E_{G,\text{noncyclic h.elem} }$. Applying 6.4.2 we see that the complement of $E_{G,\t cyclic }$ in the realization of $E_{G,\t h.elem }$ deformation retracts (in fact $G$-equivariantly) to $E_{G,\text{noncyclic h.elem} }$.

\subsubhead 6.5 Sample application\endsubsubhead 
Define ``trans-cyclic $K$-theory'' $\Bbb K^{tc}(G;F,R)$  to be the cofiber of the map induced on homology by the inclusion of cyclic isotropy into hyperelementary isotropy:
$$\Bbb H(E_{G,\text{cyclic}}/G;\Bbb K(p,R))\to \Bbb H(E_{G,\text{h.elem}}/G;\Bbb K(p,R))\to \Bbb K^{tc}(G;F,R).$$
 The iterated-homology argument  used in 5.2 for the inclusion of finite into hyperelementary expresses this  as homology: $\Bbb H(E_{G,\t h.elem }/G;\Bbb K^{tc}(p;R))$, where $\Bbb K^{tc}(p;R)$ is the coefficient system with fiber $\Bbb K^{tc}(H;F,R)$ over a point with isotropy $H$. The coefficients vanish over the subspace $E_{G,\t cyclic }/G$. A subspace over which the coefficients vanish can be deleted without changing the homology. However in 6.4.3 we saw that the join structure in the universal spaces gives an equivariant deformation retraction of the complement to the universal space $E_{G,\text{noncyclic h.elem} }$. A little more thought shows the deformation is covered by a morphism of the coefficient system, so induces an equivalence of homology $\Bbb H(E_{G,\t h.elem }/G;\Bbb K^{tc}(p;R))\simeq \Bbb H(E_{G,\text{noncyclic h.elem} }/G;\Bbb K^{tc}(p;R))$. Finally recall that a group that is both $p$ and $q$ hyperelementary for $p,q$ distinct primes must be cyclic. Thus the class of non-cyclic hyperelementary groups is the disjoint union over $p$ of non-cyclic $p$-hyperelementary groups.  All together we have shown:
\proclaim{6.5.1 Corollary} $$\Bbb K^{tc}(G;F,R)=\bigoplus_p\Bbb H(E_{G,\text{noncyclic $p$-h.elem}}/G;\Bbb K^{tc}(p;R))$$
\endproclaim
Proposition 3.2.1 implies that the $p$ summand is $p$-torsion.

 \bye